\documentclass[12pt,leqno]{amsart}
\usepackage{graphicx}
\usepackage{amssymb}
\usepackage{color}

\numberwithin{equation}{section}
\newdimen\vintkern\vintkern12pt
\def\vint{-\kern-\vintkern\int}

\newtheorem{thm}{Theorem}[section]

    \newtheorem{lem}[thm]{Lemma}
\newtheorem{cor}[thm]{Corollary}
\newtheorem{prop}[thm]{Proposition}

\newtheorem{quest}{PROBLEM}

\newtheorem{introthm}{Theorem}

 \theoremstyle{definition}
\newtheorem{definition}[thm]{Definition}
\theoremstyle{remark}
\newtheorem{rem}{Remark}[section]

\newcommand{\tref}[1]{Theorem~\ref{#1}}
\newcommand{\cref}[1]{Corollary~\ref{#1}}

\newcommand{\R}{\mathbb{R}}
\newcommand{\N}{\mathbb{N}}
\newcommand{\hm}{{\mathcal H}}
\newcommand{\Area}{\operatorname{Area}}
\newcommand{\ap}{\operatorname{ap}}
\newcommand{\Id}{\mathrm{Id}}
\newcommand{\md}{\operatorname{md}}
\newcommand{\apmd}{\ap\md}

\newcommand{\trace}{\operatorname{tr}}

\def\a{\operatorname{Area}}

\def\Jac{\operatorname{Jac}}

\def\tr{\operatorname{tr}}
\def\Fill{\operatorname{Fill}}

\def\D{\partial}

\def\R{\mathbb R}

\def\N{\mathbb N}
\def\Ha{\mathcal H}

\def\al{\alpha}

\def\ka{\kappa}

\def\eps{\epsilon}
\def\ga{\gamma}
\def\Ga{\Gamma}

\def\la{\lambda}
\def\La{\Lambda}

\def\om{\omega}
\def\Om{\Omega}

\def\si{\sigma}

\newcommand{\bt}{\begin{thm}}
\newcommand{\et}{\end{thm}}

\newcommand{\deltalip}{\delta^{\operatorname{Lip}}}

\def\wlim{\mathop{\hbox{$\om$-lim}}}

\def\<{\langle}
\def\>{\rangle}

\newcommand{\bcl}{\begin{claim}}
\newcommand{\ecl}{\end{claim}}
\newcommand{\bcor}{\begin{cor}}
\newcommand{\ecor}{\end{cor}}
\newcommand{\bdfn}{\begin{definition}}
\newcommand{\ben}{\begin{enumerate}}
\newcommand{\bit}{\begin{itemize}}
\newcommand{\blem}{\begin{lem}}
\newcommand{\bslem}{\begin{slem}}
\newcommand{\bprop}{\begin{prop}}
\newcommand{\bthm}{\begin{thm}}
\newcommand{\edfn}{\end{definition}}
\newcommand{\een}{\end{enumerate}}
\newcommand{\eit}{\end{itemize}}
\newcommand{\elem}{\end{lem}}
\newcommand{\eslem}{\end{slem}}
\newcommand{\eprop}{\end{prop}}
\newcommand{\ethm}{\end{thm}}

\usepackage{xcolor}

\begin{document}
\pagebreak


\title{Isoperimetric inequalities vs. upper   curvature bounds}

\author{Stephan Stadler}

\address
  {Max Planck Institut f\"ur Mathematik\\ 
Vivatsgasse 7\\
53111 Bonn,
Germany}
\email{stadler@mpim-bonn.mpg.de}

\author{Stefan Wenger}

\address
  {Department of Mathematics\\ University of Fribourg\\ Chemin du Mus\'ee 23\\ 1700 Fribourg, Switzerland}
\email{stefan.wenger@unifr.ch}

\date{\today}

\begin{abstract}
The Dehn function 
of a metric space measures the area necessary
in order to fill a closed curve of controlled length by a disc.
As a main result, we  prove that a length space has  curvature bounded above by $\ka$ in the sense of Alexandrov if and only if  its Dehn function is bounded above by the Dehn function 
of the model surface of constant curvature $\ka$. This extends work of Lytchak and the second author \cite{LWcurv}  from locally compact spaces to the general case. A key ingredient in the proof is the construction of minimal discs
with suitable properties in certain ultralimits. 
Our arguments also yield quantitative local and stable versions of our main result. 
The latter has implications on the geometry of asymptotic cones.
 \end{abstract}

\maketitle

\renewcommand{\theequation}{\arabic{section}.\arabic{equation}}
\pagenumbering{arabic}

\section{Introduction}
\subsection{Main result}
The main result of this paper is an analytic characterization of CAT($\kappa$) spaces -- complete metric spaces with curvature bounded above
by $\kappa$ in the sense of Alexandrov. For simplicity, we first describe the result for CAT(0) spaces.
 We say that a metric space $X$ satisfies the Euclidean isoperimetric inequality for curves if any closed Lipschitz curve $\gamma$ in $X$  bounds  a Lipschitz 
disc $v$ in $X$ such that
\[\a(v)\leq\frac 1 {4\pi}\cdot \ell ^2(\gamma ). \]
Here, $\a(v)$ denotes the parametrized Hausdorff area of $v$ and $\ell(\gamma )$ is the length of $\gamma$.
The fact that CAT(0) spaces satisfy the Euclidean isoperimetric inequality for curves
is a well-known consequence of Re\-shetnyak's majorization theorem \cite{Resh-majorization} and the isoperimetric inequality in the Euclidean plane.
Vice versa, for locally compact metric spaces $X$, it was proved in \cite{LWcurv} that a Euclidean isoperimetric inequality implies the zero upper curvature
bound. Moreover, in \cite{LWcurv}, this isoperimetric characterization of upper curvature bounds was extended to non-zero bounds.
In order to formulate our main result, it is convenient to introduce 
{\em Dehn functions}. 
In general, if $X$ is a metric space, a Dehn function $\delta_X(r)$ controls the area
needed in order to fill a closed curve of length at most $r$ in $X$ by a disc.
A precise definition requires a choice of area and types of discs.

For instance, chosing parametrized Hausdorff area and Lipschitz discs results in the {\em Lipschitz Dehn function} $\deltalip_X$.
More precisely, 
let $\delta:(0,\infty)\to[0,\infty]$ be a function such that any Lipschitz circle of length at most $r$ bounds a Lipschitz disc of area at most $\delta(r)$. Then $\deltalip_X$ is the greatest lower bound for all such functions.

Denote by $M^2_\kappa$ the complete simply connected surface of constant sectional curvature $\kappa$. Thus, up to scaling, we see the round 2-sphere
$M^2_1=\mathbb{S}^2$, the Euclidean plane $M^2_0=\R^2$, and the hyperbolic plane $M^2_{-1}=\mathbb{H}^2$. 
Let us denote by $D_\kappa$ the diameter of $M^2_\kappa$ and by $\delta_\kappa$ 
the Lipschitz Dehn function of $M^2_\kappa$.

\begin{introthm}\label{thm_main}
Let $X$ be a complete length space and $r_0\in(0,2D_\kappa]$. Suppose that
\begin{equation}\label{eq_isop}
    \deltalip_X(r)\leq\delta_\ka(r)
\end{equation}
holds for all $r\in(0,r_0)$. Then every closed ball of radius at most $\frac{r_0}{4}$
in $X$
is convex and CAT($\ka$). Moreover, if $r_0=2D_\ka$, then $X$ is a CAT($\ka$)
space.
\end{introthm}

This extends the main theorem in \cite{LWcurv} from locally compact spaces to the general case and at the same time provides a quantitative local version. For even more general results, see Section~\ref{sec_general} below.
Accordingly, CAT($\kappa$) geometry is as much an analytic theory as it is a geometric theory.
In particular, upper curvature bounds can be identified without the knowledge of geodesics or angles.
One encounters such situations in many geometric settings \cite{AB04,AB16,LS_conf,LS_imp,LS_curv,LWa_curv,PS,RR_2D}. 

Beckenbach and Rado in \cite{BR_sub} first discovered a relationship between isoperimetric inequalities and upper curvature bounds on smooth 2-dimensional Riemannian manifolds.
Later, this was generalized to non-smooth metric surfaces by Reshetnyak \cite{R_iso}.
In \cite{LWcurv}, for locally compact spaces, the isoperimetric inequality is translated to an upper curvature bound
in
three steps: 1) Solve the Plateau problem for a given Jordan curve; 2) Show that the solution is a minimal disc satisfying the same isoperimetric inequality as the surrounding space; 
3) Prove that the minimal disc is intrinsically a metric surface with the desired upper curvature bound.
\medskip

If we were able to solve the Plateau problem in the setting of Theorem~\ref{thm_main},
then the strategy of \cite{LWcurv} would still be successful. Inspection of the proof shows
that it would even be enough to solve the following problem: 
For a given Jordan curve $\Ga\subset X$ construct 
a larger space $Y\supset X$ which satisfies \eqref{eq_isop}
and such that the Plateau problem for $\Ga$ is solvable in $Y$.

Even though there are natural candidates for the space $Y$,
such as ultra-completions $X_\om$ of $X$,  we are unable to do this. As a matter of fact, it is possible to construct
for a given Jordan curve $\Ga\subset X$ a minimal disc in $X_\om$ filling $\Ga$ but we do not know whether $X_\om$ satisfies inequality \eqref{eq_isop}, see Section~\ref{sec:questions}.
Despite these obstacless, our proof of Theorem~\ref{thm_main} still employs minimal surface theory.
We sidestep the fact that we are unable to answer this question
by producing minimal discs in ultralimits whose intrinsic structure appears as if the Dehn function of the surrounding space would satisfy the correct bounds, see the next section.

\subsection{Minimal surfaces in ultralimits of locally non-compact spaces}
As announced, we prove a general result which serves as an appropriate substitute for the solvability of the Plateau Problem in locally compact spaces. 
This is the technical heart behind Theorem~\ref{thm_main} and its generalizations discussed below.
We state it here in a very simplified setting, for the version in full generality see Theorem~\ref{thm_plateau}. In the context of the Plateau problem, it is natural to enlarge the class of admissible discs filling a given circle from Lipschitz
to Sobolev discs, see Section~\ref{sec_sob_basics} for precise definitions.
As in the classical case of Riemannian manifolds, this leads to better compactness properties relative to energy bounds.
From now on we will focus on the (Sobolev) Dehn function $\delta_X(r)$ instead of the Lipschitz Dehn function $\deltalip_X(r)$. Since every Lipschitz disc is Sobolev, we obtain the natural inequality $\delta_X(r)\leq\deltalip_X(r)$. Thus, results which only involve upper bounds on  the Dehn function are a priori stronger than their Lipschitz counterparts.

\begin{introthm} \label{thm_minimal}
Let $X$ be a complete length space, $r_0>0$ and $\kappa\in\R$. Suppose $\Ga\subset X$ is a Jordan curve 
with $\ell(\Ga)<r_0$ and the Dehn function of $X$ satisfies $\delta_X(r)\leq\delta_\kappa(r)$
for all $r\in(0,r_0)$. Then there exists an ultracompletion $X_\om$ of $X$ and a
continuous map $v:\bar D\to X_\om$ which is a 
solution to the Plateau problem for $\Ga$ in $X_\om$ and satisfies
\begin{equation}\label{eq_minimal}
\a(v|_\Om)\leq\delta_\ka(\ell(v|_{\D\Om}))
\end{equation}
for every Jordan domain $\Om\subset D$ such that $\ell(v|_{\D\Om})<r_0$.
\end{introthm}

Recall that an ultracompletion $X_\om$ is a certain metric space that contains an isometric  copy of $X$ and which is constructed with the help of a non-principal ultrafilter $\om$ on the natural numbers, see Section~\ref{sec_ultra} for details.  Note that in Theorem~\ref{thm_minimal}
we do not gain control on the Dehn function $\delta_{X_\om}$ of $X_\om$. 
However, the ``intrinsic isoperimetric inequality'' \eqref{eq_minimal} 
has the effect that $v$ behaves as if $\delta_{X_\om}$ would be bounded above by $\delta_\ka$.

Here is how we deduce Theorem~\ref{thm_main} from Theorem~\ref{thm_minimal} in the case that $X$ is geodesic and $\ka=0$:
For a Jordan triangle $\triangle\subset X$ we obtain from Theorem~\ref{thm_minimal}  a minimal disc $v$ filling $\triangle$
in some ultracompletion $X_\om$.
As follows from \cite{LWint}, the map $v$ factors as
\[D\xrightarrow{P} Z_v\xrightarrow{\bar v}X_\om   \]
where $Z_v$ is a metric disc and $\bar v$ is a 1-Lipschitz map which restricts to 
an arclength preserving homeomorphism $\D Z_v\to\Ga$. 
Moreover,  for every Jordan domain $O\subset Z_v$ holds
\[\hm^2(O)\leq\delta_\ka(\ell(\D O)).\]
We deduce from \cite{LWcurv} that $Z_v$ is CAT(0) and then
Reshetnyak's majorization theorem implies that $X$ itself is CAT(0).

\subsection{Generalizations}\label{sec_general}
We prove a stable version of Theorem~\ref{thm_main} for sequences
of metric spaces. The setting involves metric spaces without any control on small scales.
This requires us to adjust the way we measure area in the definition of the Dehn function. Instead of parametrized Hausdorff area we will use the {\em Riemannian inscribed area} $\a_{\mu^i}$, originally defined by Ivanov in \cite{I_vol}, and the
associated {\em Riemannian Dehn function} $\delta_X^{\mu^i}$.
Informally speaking, instead of assigning to a unit ball in a normed plane its
area, one assigns the maximal area of an inscribed ellipsoid -- the {\em John ellipsoid}.
The precise definition and a discussion of basic properties can be found in Section~\ref{sec_length}. Here we only mention that for many geometrically interesting spaces
these subtleties disappear and the equality $\a_{\mu^i}=\a$ holds.
For instance, this is the case for all metric spaces with curvature locally bounded either above or
below, in particular, it holds for all Riemannian manifolds.

\begin{introthm}\label{thm_intro_stable}
Let $(X_n)$ be a sequence of complete length spaces and $\kappa\in\R$.
Suppose $r_0\in(0,2D_\kappa]$ and
\[\delta^{\mu^i}_{X_n}(r)\leq(1+\eps_n)\cdot\delta_\ka(r)+\eps_n\]
holds for all $r\in(0,r_0)$ and some sequence $\eps_n\to 0$. Then every ultralimit $X_\om$ is locally CAT($\ka$). More precisely, every closed ball of radius at most
$\frac{r_0}{4}$ in $X_\om$ is convex and CAT($\ka$).
\end{introthm}

Note that Theorem~\ref{thm_main} really is a special case.
Indeed, in the setting of Theorem~\ref{thm_main} we have $\deltalip_X(r)=\delta_X(r)=\delta_X^{\mu^i}(r)$. Then Theorem~\ref{thm_intro_stable} applies to the constant sequence
$X$ and provides the curvature bound for all of its ultracompletions and therefore for $X$
itself.

The result also has implications for asymptotic cones. 

\begin{introthm}\label{thm_asymp_cones}
Let $X$ be a complete length space  
such that
\[\limsup\limits_{r\to\infty}\frac{\delta_X^{\mu^i}(r)}{r^2}\leq\frac{1}{4\pi}.\]
Then every asymptotic cone of $X$ is a CAT(0) space. Moreover, if the inequality is strict, then every asymptotic cone of $X$ is a tree. In particular, in this case, $X$
is Gromov hyperbolic.
\end{introthm}

The first statement is a version of \cite[Theorem~1.1]{W_asym} for locally non-compact spaces and the second statement is a variant of \cite[Theorem~1.1]{W_gromov}.
For a localized version of the second statement in  Theorem~\ref{thm_asymp_cones} 
see Theorem~\ref{thm_tree}.

Theorem~\ref{thm_intro_stable} can also be used to turn fine infinitesimal information on
the Dehn function into curvature bounds:

\begin{introthm}\label{thm_fine}
Let $\delta:(0,r_0)\to \R$ be a continuous non-decreasing function with
\[\limsup\limits_{r\to 0}\frac{\delta(r)-\frac{r^2}{4\pi}}{r^4}\leq 0.\] 
Suppose that  $(X_n)$ is a sequence of complete length spaces
such that the Riemannian Dehn functions satisfy 
\[\delta^{\mu^i}_{X_n}(r)\leq\delta(r)+\eps_n\] on 
$(0,r_0)$ for some sequence $\eps_n\to 0$. Then any ultralimit
$X_\om$ is locally CAT(0). More precisely, there exists $\tilde r>0$ depending only on 
the function $\delta$ such that every closed ball in $X_\om$ of radius at most $\tilde r$ is
convex and CAT(0).
\end{introthm}

\subsection{Motivation and strategy for Theorem~\ref{thm_minimal}}
For simplicity, we restrict this discussion to the case $\kappa=0$.

Recall that in the non-locally-compact case we are unable to solve the Plateau problem
in the traditional sense.
However, if $X$ was  known to be CAT(0), then the Plateau problem would be solvable, $X$ being locally compact or not. For this one chooses a minimizing sequence of discs filling a given Jordan curve $\Ga\subset X$ and constructs a suitable limit $v$ in an ultracompletion $X_\om$. Since $X$ is CAT(0) there exists a 1-Lipschitz retraction $X_\om\to X$ and hence we can push $v$ back to $X$ without increasing energy or area \cite{GW_noncompact,St}.

While for a locally non-compact space $X$ which satisfies the Euclidean isoperimetric inequality we are not able to a priori show the existence of such retractions,
this still motivates the search for minimal discs in ultracompletions.

To prove Theorem~\ref{thm_minimal} we start with a minimizing sequence $(v_n)$ of Sobolev discs
filling $\Ga\subset X$. 
We then reparametrize using Morrey's $\eps$-conformality lemma to
make energy and area almost equal.
Using Rellich-Kondrachov compactness, we  select a subsequence which $L^2$-converges
to a Sobolev disc $v$
in some auxiliary metric space.
Then comes the critical part:
We use the isoperimetric information to show that the filling area of $\Ga$ cannot drop
in $X_\om$. Using our setup, this allows us to show that the limit  is in fact a minimal disc and that areas and energies converge.
 Once this is established, we obtain strong convergence in the Sobolev norm.
Next we apply Fuglede's lemma to see length convergence $\ell(v_n|_\ga)\to\ell(v|_\ga)$ for most paths $\ga\subset D$.
Using area convergence and the bound on the Dehn function we conclude
the intrinsic isoperimetric inequality for the limit $v$, thus Theorem~\ref{thm_minimal}.

\subsection{Further questions}\label{sec:questions}
As already mentioned above, it is not known whether an isoperimetric inequality passes to ultracompletions. More generally, we pose the following problem:

\begin{quest}\label{quest_stable}
Let $\delta:(0,r_0)\to\R$ be a continuous non-decreasing function.
Let $(X_n)$ be a sequence of complete length spaces
such that the Dehn functions satisfy
\[\delta_{X_n}(r)\leq\delta(r)\]
on $(0,r_0)$. Does the Dehn function of an
ultralimit of $(X_n)$ satisfy the same inequality?
\end{quest}

On a more technical level, it is also not clear whether the use of the Riemannian inscribed area in Theorem~\ref{thm_intro_stable} is really necessary. The crucial point in the proof where we use it is in the form of Morrey's $\eps$-conformality lemma which allows us to reparametrize a Sobolev disc such that 
area and energy become almost equal. Morrey's lemma fails for the Hausdoff area as can be seen in non-Euclidean normed spaces. However, constructing a counterexample to 
Theorem~\ref{thm_main} where the Riemannian Dehn function is replaces by the Hausdorff
Dehn function seems very difficult. The examples would have to involve locally non-compact  spaces.

\subsection{Organization}
In Section~\ref{sec_pre} we set notation and collect the necessary background from metric geometry. In Section~\ref{secSob} we recall relevant parts of the Sobolev theory in metric spaces. In several subsections we discuss
 energy, area, isoperimetric inequalities, quasi-conformality and regularity results.
 In Section~\ref{sec_dehn} we introduce Dehn functions, construct universal thickenings of metric spaces and recall the Plateau problem and the intrinsic structure of its solutions.
 Section~\ref{sec:ET} is devoted to an infinitesimal Euclidean property for metric spaces -- property (ET) -- and how to ensure it holds, given some isoperimetric control. In Section~\ref{sec_fill} we begin working towards our main theorem.
 We prove a result (Proposition~\ref{prop_lift_fillingarea}) which ensures that the filling area cannot drop in an ultralimit
 of a sequence of metric spaces, assuming certain bounds on the Dehn functions.
 In Section~\ref{sec_int} we prove a technical result (Proposition~\ref{prop_inisin})
 stating that a certain limit of Sobolev discs satisfies an intrinsic isoperimetric inequality
 if each individual Sobolev disc lives in a space which supports an isoperimetric inequality.
 In Section~\ref{sec_plateau} we solve a version of the Plateau problem for locally non-compact spaces (Theorem~\ref{thm_plateau}). In the final Section~\ref{sec_main} we provide proofs of our main results.

\subsection{Acknowledgements}
We like to thank Alexander Lytchak for several inspiring discussions.
The first author was supported by DFG grant SPP 2026. The second author was partially supported by Swiss National Science Foundation grant 212867.

\section{Basics on metric spaces} \label{sec_pre}
\subsection{Notation}
 Connected open subsets of $\R^n$ will be called {\em domains}.
A domain $\Omega\subset\R^n$ is  a {\em Lipschitz domain}, if its boundary $\D\Om$ can locally be written as the graph of a Lipschitz function. When $n=2$ then this is equivalent to the requirement that $\partial \Omega$ is locally bilipschitz to an open interval \cite{Tuk80}. The open unit disc in $\R^2$ is denoted by $D$ and the standard annulus $S^1\times[0,1]$ by $A$. We call annulus any domain $U\subset\R^2$ homeomorphic to $A$.
We denote the two {\em boundary circles} of an annulus $U\subset \R^2$ by $\D^\pm U$. By \cite{Tuk80}, $U$ is bilipschitz to $A$ if $\D^\pm U$ are bilipschitz curves.

We will denote distances in a metric space $X$ by $d$ or $d_X$.
Let $X=(X,d)$ be a metric space. The open ball  in $X$ of radius $r$ and center $x_0\in X$ is denoted by
\[B_r(x_0)= \{x\in X: d(x_0, x)<r\}. \]
More generally, for any subset $A\subset X$ we denote its open tubular neighborhood of radius $r$ by
\[N_r(A)=\{x\in X: d(A, x)<r\}.\]

A {\em Jordan curve} in $X$ is a subset $\Gamma\subset X$ which is homeomorphic to $S^1$. Given a Jordan curve $\Gamma\subset X$, a continuous map $\gamma: S^1\to X$
with image $\Ga$ is called a {\em weakly monotone parametrization} of $\Gamma$ if it has connected fibers.

For $m\geq 0$, the $m$-dimensional Hausdorff measure on $X$ is denoted by $\mathcal H^m=\mathcal H^m_X$. The normalizing constant is chosen in such a way that on Euclidean space $\R^m$ the Hausdorff measure $\mathcal H^m$ equals the Lebesgue measure.

The length of a curve $\gamma $ in a metric space $X$ will be denoted by $\ell_X(\gamma )$ or simply by $\ell(\gamma)$.
 A continuous curve of finite length is called \emph{rectifiable}.
A (local) \emph{geodesic} in a space $X$ is a (locally) isometric map  from an interval to $X$.
   A space $X$ is called \emph{a geodesic space} if any pair of points in $X$ is connected by a geodesic.
	
	For $\eps>0$ we call a Lipschitz curve $c:[a,b]\to X$  an {\em $\eps$-geodesic},
if it has constant speed and satisfies 
\[\ell(c)\leq(1+\eps)\cdot d(c(a),c(b)).\]
	
A space $X$ is a \emph{length space}
if for all $x,y\in X$ the distance $d(x,y)$ equals $\inf \{\ell_X (\gamma ) \}$, where
$\gamma$ runs over the set of all curves connecting $x$ and $y$.
In a length space any pair of points is connected by an $\eps$-geodesic for every $\eps>0$.

\subsection{Intrinsic metric of a map} \label{sec_associated}
We refer the reader to \cite{BBI, LWint,PS} for  discussions of the following construction and related topics.
Let $Z$ be a topological space and $X$ a metric space. Let $u:Z\to X$ be a continuous map. The {\em intrinsic distance associated with $u$}
is the function $d_u:Z\times Z\to [0,\infty ]$
 defined by
\[d_u(z_1,z_2)=\inf\{\ell_X(\ga)|\ \ga \text{ is a path in } Z \text{ connecting }z_1\text{ and }z_2\}.\]
If it only takes finite values, then it defines a pseudo-metric. The associated metric space $Z_u$
which arises from  identifying pairs of points at zero $d_u$-distance, is a length
space.  We will call it the \emph{intrinsic metric space associated with the map $u$}.

By construction, the space $Z_u$ associated with  the map $u$ comes with a canonical, possibly non-continuous, surjective projection $P:Z\to Z_u$ and
a $1$-Lipschitz map $\bar u:Z_u\to X$ such that $u=\bar u\circ P$.

If $X$ is a metric space in which any pair of points is connected by a curve of finite length,
then the \emph{length space $X^i$ associated to $X$} is the special case $X^i=Z_u$ of the above construction where $u$ is the identity map $u=\Id:X\to X$.
The completeness of $X$ implies that $X^i$ is complete as well. The $1$-Lipschitz map $\bar u :X^i\to X$ from above  is the identity in this case. The map $P =\bar u ^{-1}:X\to X^i$  need not
 be continuous, but it sends curves of finite length in $X$ to continuous curves of the same length in $X^i$.

\subsection{Ultralimits}\label{sec_ultra}
We refer the reader to \cite{AKP_found} for an extended treatment of ultralimits in the context of metric geometry. 
For a non-principal ultrafilter $\om$ on $\N$ and a sequence of pointed  metric spaces $(X_n, x_n)$, we will often consider their ultralimit
\[(X_\om,x_\om)=\wlim(X_n,x_n)\]
which is a pointed metric space whose elements are equivalence classes of sequences $(p_n)$ which are bounded relative to $(x_n)$. The metric $d_\om$ on $(X_\om,x_\om)$
is induced by the metrics $d_n$ on $X_n$ via $d_\om((p_n),(q_n))=\wlim d_n(p_n,q_n)$.
Often the choice of basepoints $x_n$ is irrelevant in our considerations and then we neglect it in our notation.
The ultralimit of the constant sequence $X$ with fixed base point $x\in X$ will be called {\em ultrapower} or {\em ultracompletion}  of $X$ with respect to $\om$. Ultrapowers do not depend on the choice of basepoint. Note that any ultralimit is a complete metric space. Moreover, any metric space admits a canonical isometric embedding into its ultrapowers via constant sequences. 
In case every $X_n$ is a length space, then any ultralimit $X_\om$ is geodesic.

For any sequence of subsets $A_n\subset X_n$ we denote by $A_\om\subset X_\om$
their ultralimit. So $A_\om\subset X_\om$ corresponds precisely to those points in $X_\om$
which can be represented as a sequence $(a_n)$ with $a_n\in X_n$ and $\sup d(a_n,x_n)<\infty$. Note that $A_\om$ is always closed, even if the individual $A_n$ might  not be. 

If $(X_n)$ and $(Y_n)$ are  sequences of metric spaces such that  $X_n\subset Y_n$ 
and $Y_n\subset N_{\eps_n}(X_n)$ for some sequence $\eps_n\to 0$,
then $\wlim(X_n,x_n)=\wlim(Y_n,x_n)$ for any choice of basepoints $x_n\in X_n$.

\blem\label{lem_Jordan_lift}
Let $X_\om$ be an ultralimit of a sequence of length spaces $X_n$.
Suppose that $\Ga_\om\subset X_\om$ is a rectifiable Jordan curve.
Then, for every $\eps>0$ there exist rectifiable Jordan curves $\Ga_n\subset X_n$ with
$\wlim \Ga_n=\Ga_\om$ and $\ell(\Ga_n)\leq(1+\eps)\cdot\ell(\Ga_\om)$ for all
$n\in\N$. 
\elem

\proof
Let $c_\om:S^1\to X_\om$ be a constant speed parametrization of $\Ga_\om$.
By \cite[Corollary~2.6]{LWY_dehn}, there is a sequence of uniformly Lipschitz maps $c_n:S^1\to X_n$ which
ultra-converges to $c_\om$, meaning $\wlim c_n(t)=c(t)$ for all $t\in S^1$, and such that $\ell(c_n)\leq(1+\eps)\cdot\ell(c_\om)$ holds for all  $n\in\N$.

Now we will find the 
desired Jordan curves $\Ga_n$ inside the images of the $c_n$.
Choose three equidistant points $t^j\in S^1, j=1,2,3$
and set $z_n^j=c_n(t^j)$. Denote by $\al^j\subset S^1$
the closure of the component of $S^1\setminus\{t^1,t^2,t^3\}$ which does not contain the point $t^j$. Next, inside the image $c_n(\al^j)$, choose a minimizing geodesic
$\ga^j_n$ from 
$z_n^{j+1}$ to $z_n^{j+2}$. Note that these geodesics ultra-converge to  $c_\om(\al^j)$.
For every $\delta>0$ and large enough $n$, the geodesics $\ga_n^j$ are disjoint away from the $\delta$-neighborhood
of $\{z_n^1,z_n^2,z_n^3\}$. Now we orient the three geodesics and pass to subsegments
$\hat\ga_n^j\subset\ga_n^j$ where the starting point of $\hat\ga_n^j$ is the last point of $\ga_n^j$ on $\hat\ga_n^{j-1}$. Finally, put $\Ga_n=\bigcup_{j=1}^3 \hat\ga_n^j$.
\qed

\subsection{CAT($\kappa$)}

Let $\ka$ be a real number. Recall that a CAT($\kappa$) space is a complete metric space $X$ where any pair of points at distance strictly less than $D_\kappa$ is joined by a geodesic and such that distances between points on a geodesic triangle $\triangle\subset X$ of perimeter strictly less than $2 D_\kappa$
are bounded above by the distances between corresponding points on the comparison triangle $\tilde\triangle\subset M^2_\kappa$.
In order to check if a geodesic space is CAT($\ka$)
it suffices to prove the CAT($\ka$) comparison for {\em Jordan triangles} --
geodesic triangles which are Jordan curves, cf. \cite[Lemma~3.1]{LWcurv}.
A more flexible characterization of CAT($\kappa$) spaces which does not refer to geodesics can be provided using majorizations.
Recall that a rectifiable Jordan curve $\Ga$ in a metric space $X$ admits a {\em $\kappa$-majorization in $X$}, if there is a closed convex region $C\subset M^2_\kappa$
and a 1-Lipschitz map $C\to X$ which restricts to an arclength preserving homeomorphism $\D C\to\Ga$. If $\ka$ is clear from the context we will simply speak of majorizations.

By Reshetnyak's majorization theorem \cite{Resh-majorization},
every rectifiable Jordan curve in a  CAT($\kappa$) space  can be  $\kappa$-majorized.
On the other hand, if a Jordan triangle in a metric space admits a $\ka$-majorization,
then it clearly satisfies the CAT($\kappa$)  comparison.
We will make use of the following, which is a consequence of \cite[Radial lemma~9.52]{AKP_found}.

\blem\label{lem_loc_cat}
Let $X$ be a complete geodesic space. Suppose there is $r_0\in(0,2D_\kappa]$ such that every Jordan triangle of perimeter 
strictly less than $r_0$ satisfies the CAT($\ka$) comparison. Then 
every closed $r$-ball with $r\leq \frac{r_0}{4}$ is convex and CAT($\ka$).
Moreover, if $r_0=2D_\ka$, then $X$ is a CAT($\ka$) space.
\elem

\section{Sobolev theory} \label{secSob}

\subsection{Basics}\label{sec_sob_basics}
We recall basic definitions of
Sobolev maps with values in  a metric space and refer to \cite{HKST15,KS93,LW_plateau,LWint,Res97} and references
therein for further information.    Let $\Omega$ be a bounded domain in $\R^n$ and $(X,d)$  a complete metric space.
For $p>1$ let $L^p(\Om, X)$ be the set of measurable and essentially separably valued maps $u: \Om\to X$ such that for some and thus every 
$x\in X$ the function $u_x(z):= d(x, u(z))$ belongs to the classical space $L^p(\Om)$ of $p$-integrable functions on $\Om$. 

\bdfn
 A map $u\in L^p(\Om, X)$ belongs to the Sobolev space $W^{1,p}(\Om, X)$ if there exists $h\in L^p(\Om)$ such that $u_x$ is in the classical Sobolev space 
$W^{1,p}(\Om)$ for every $x\in X$ and its weak gradient satisfies $|\nabla u_x|\leq h$ almost everywhere.
The localized spaces $W^{1,p}_{loc}(\Om)$ are defined similarly.
\edfn

There are several natural notions of energy of a map $u\in W^{1,p}(\Om, X)$. Throughout this paper we will use the {\em Reshetnyak $p$-energy} defined by 
\[E_+^p(u):= \inf\left\{\|h\|_{L^p(\Om)}^p\;\big|\; \text{$h$ as in the definition above}\right\}.\]

Recall that a family $\mathcal{C}$ of curves in $\Om$ is called {\em $p$-exceptional},
if
there exists a $p$-integrable Borel function $\si:\Om\to[0,\infty]$
such that for every locally rectifiable  curve $\ga\in\mathcal C$ the path integral satisfies 
\[\int_\ga \si\, ds=\infty\, .\]
We say that a property holds for {\em $p$-a.e.~curve} in $\Om$ if
the family of curves on which the property fails is $p$-exceptional.
  For instance, if a property holds for  $p$-a.e.~curve in $\Om$ and  $F:[0,1]^{n-1}\times[0,1]\to \Omega$ is a bilipschitz embedding, then  for almost all 
$x\in [0,1]^{n-1}$  the property holds true for the curve 
$\gamma _x(t)=F(x,t)$.

 A map $u\in L^p (\Omega,X)$ is contained in the Sobolev space $W^{1,p} (\Omega,X)$
if and only if there exist a Lebesgue representative $\bar u$ of $u$ and a Borel  function $\rho \in L^p (\Omega )$
such that for \emph{$p$-a.e.~curve} $\gamma: [0,1] \to \Omega$  the composition $\bar u\circ \gamma$ is  continuous  and
\begin{equation} \label{eq-n1p}
 \ell_X( \bar u\circ \gamma ) \leq \int _{\gamma} \rho\, ds.
\end{equation}
In what follows, we will always choose such a representative $\bar{u}$ of $u$ and will simply denote it by $u$.

 There exists a  minimal function  $\rho=\rho _u$ satisfying the condition above, uniquely defined up to sets of measure zero. It will be called  the 
\emph{generalized  gradient} or {\em minimal weak upper gradient} of $u$.  By \cite{LW_plateau,Res97}, the $p$-th power of the $L^p$-norm of $\rho_u$
 coincides  with the  Reshetnyak $p$-energy defined above, 
	\[E_+ ^p (u)=\|\rho_u\|^p_{L^p(\Om)}=\int _{\Omega} \rho ^p _u (z)\, dz.\]

If $\Om$ is a Lipschitz domain, then every $u\in W^{1,p}(\Om, X)$ has a canonically defined {\em trace} $\tr(u)\in L^p(\D\Om,X)$, cf. \cite{KS}. For instance, if $\Om$
is the open unit ball in $\R^n$, then for almost every $z\in S^{n-1}$ the map $t\mapsto u(tz)$ is in $W^{1,p}((\frac{1}{2},1),X)$ and
\[\tr(u)(z)=\lim\limits_{t\to 1}u(tz).\]
For a general Lipschitz domain $\Om$, if $u$ has a continuous extension $\hat{u}$ to $\overline{\Om}$ then $\trace(u)$ is just the restriction of $\hat{u}$ to $\D\Om$. 

If a Lipschitz domain $\Om$ is a union of two disjoint Lipschitz subdomains $\Om^\pm$ and the Lipschitz boundary $T=\D\Om^-\cap\D\Om^+$, and  $u^\pm\in W^{1,2}(\Om^\pm,X)$ have the same trace on $T$, then one obtains a  map $u\in W^{1,2}(\Om,X)$
by gluing $u^\pm$ along $T$, see \cite[Theorem~1.12.3]{KS}.

\subsection{Length, energy, and area}\label{sec_length}

Every map $u\in W^{1,p}(\Om, X)$ has an {\em approximate metric derivative} at almost every point $z\in \Om$ in the following sense, see \cite{Kar_sobolev} and \cite{LW_plateau}. There exists a unique semi-norm on $\R^n$, denoted $\apmd u_z$, such that 
 \begin{equation*}
    \ap\lim_{z'\to z}\frac{d(u(z'), u(z)) - \apmd u_z(z'-z)}{|z'-z|} = 0,
 \end{equation*}
where $\ap\lim$ denotes the approximate limit, see \cite{Evans}. If $u$ is Lipschitz, then the approximate limit can be replaced by an honest limit.
 The map $z\mapsto \apmd u_z$ into the space of semi-norms has a Borel measurable representative \cite{LW_plateau}.
For $p$-a.e.~absolutely continuous curve $\gamma:I\to \Omega$ we have:
 \begin{equation} \label{eq:almostall}
\ell_X(u\circ \gamma )=\int _I \apmd u_{\gamma (t)}  (\gamma '(t)) dt.
\end{equation}
Moreover, for almost every  $z\in \Omega$ we have $\rho _u (z) = \sup\limits_{v\in S^{n-1}} \apmd u_z (v)$. It follows from \cite{LW_plateau} that 
\[E_+^p(u) = \int_\Om\mathcal{I}_+^p(\apmd u_z)\,dz,\] 
where for a semi-norm $s$ on $\R^n$ we have set 
\[\mathcal{I}_+^p(s):= \max\{s(v)^p : |v|=1\}.\]

Later on, we will make use of the following lemma which is a slight variant of \cite[Theorem~7.3.9]{HKST15}.

\blem\label{lem_lim}
Let $(u_n)_{n\in\N}$ be a sequence in $W^{1,p}(\Om, X)$ which converges in $L^p(\Om,X)$
to a function $u$. Suppose that the corresponding sequence $(\rho_{u_n})_{n\in\N}$ of minimal weak upper gradients 
is uniformly bounded in $L^p(\Om)$ and converges weakly in $L^p(\Om)$
to a function $\rho$. Then $u\in W^{1,p}(\Om,X)$ and $\rho$ is a weak upper gradient of $u$. In particular,
\[\|\rho_u\|_{L^p(\Om)}\leq\liminf_{n\to\infty} \|\rho_{u_n}\|_{L^p(\Om)}.\]
\elem

\proof
Let us isometrically embed $X$ into a Banach space $V$.
By Mazur's lemma \cite[Section~2.3]{HKST15}, we can form sequences $(\tilde u_l)_{l\in \N}$ and $(\tilde\rho_l)_{l\in\N}$
of convex combinations of $(u_n)_{n\in\N}$ and $(\rho_{u_n})_{n\in\N}$, respectively, such that 
$\tilde u_l\to u$ in $L^p(\Om,V)$, $\tilde\rho_l\to\rho$ in $L^p(\Om)$ and for every $l\in\N$ and $p$-a.e.~rectifiable path
$\ga$ in $\Om$ holds
\[\ell_V(\tilde u_l\circ\ga)\leq\int_\ga \tilde\rho_l\, ds.\]
By \cite[Proposition~7.3.7]{HKST15}, $u\in W^{1,p}(\Om,V)$ with $\rho$
as a weak upper gradient. Since a subsequence of $(u_n)_{n\in\N}$ converges pointwise almost everywhere to $u$, we 
obtain $u\in W^{1,p}(\Om,X)$. The last statement follows from the semi-continuity of the norm with respect to weak convergence 
since $\|\rho_u\|_{L^p(\Om)}\leq\|\rho\|_{L^p(\Om)}$.
\qed

\medskip

We will mostly be interested in Sobolev discs and, more generally, Sobolev images of 
planar domains. So let us for the rest of this section specialize to the case where $\Om$ is a domain
in $\R^2$.

\bdfn
The (parameterized Hausdorff) area of a map $u\in W^{1,2}(\Om, X)$ is defined by 
\[\a(u):= \int_\Om \Jac(\apmd u_z)\,dz,\] 
where the Jacobian $\Jac(s)$ of a semi-norm $s$ on $\R^2$ is the Hausdorff $2$-measure in $(\R^2, s)$ of the Euclidean unit square if $s$ is a norm and $\Jac(s)=0$ otherwise. 
\edfn

We will furthermore need a somewhat different definition of paramet\-rized area, also known as the {\em inscribed Riemannian area} defined as follows. The $\mu^i$-Jacobian $\Jac_{\mu^i}(s)$ of a norm $s$ on $\R^2$ is given by $$\Jac_{\mu^i}(s):= \frac{\pi}{|L|},$$ where $|L|$ is the Lebesgue measure of the ellipse of maximal area contained in $\{v\in\R^2: s(v)\leq 1\}$. If $s$ is a degenerate semi-norm then we set $\Jac_{\mu^i}(s)=0$. The inscribed Riemannian area of a map $u\in W^{1,2}(\Om, X)$ is defined by $$\a_{\mu^i}(u):= \int_\Om \Jac_{\mu^i}(\apmd u_z)\,dz.$$
The two notions of area are related by 
$$\frac{\pi}{4}\cdot \a_{\mu^i}(u) \leq \a(u)\leq \a_{\mu^i}(u)$$ for every $u\in W^{1,2}(\Om, X)$, see \cite[Section 2.4]{LW_energy}. The following semi-continuity properties of energy and area are essential.

\bprop[{\cite[Theorem~4.2]{Res97}, \cite[Corollary~5.8]{LW_plateau}}]\label{prop_semi}
For $p\geq 2$, let $(u_n)_{n\in\N}$ be a sequence in $W^{1,p}(\Om, X)$ which converges in $L^p(\Om,X)$
to $u\in W^{1,p}(\Om, X)$.
Then, 
\[E_+^p(u)\leq\liminf_{n\to\infty} E_+^p(u_n) \text{ and } \a(u)\leq\liminf_{n\to\infty} \a(u_n).\]
The second inequality moreover holds with the Hausdorff area replaced by the inscribed Riemannian area.
\eprop

The following definition is of central importance.
\bdfn
A metric space $X$ supports a {\em $(C,l_0)$-isoperimetric inequality}, if
for every Lipschitz curve $\ga:S^1\to X$ of length $l< l_0$ there exists a Sobolev
disc $u\in W^{1,2}(D,X)$ with $\tr(u)=\ga$ and 
\[\a(u)\leq C\cdot l^2.\]
\edfn

Note that changing the definition above from Hausdorff to Riemannian inscribed area
will only change the constant $C$.

Isoperimetric inequalites as above have an important effect on the regularity of
minimal discs \cite{LW_plateau}. 
Moreover, such inequalites allow for flexible reparametrizations.
To explain this, recall from \cite{LW_plateau} the following terminology. Let $T\subset \R^2$
be a subset which is bilipschitz to an open interval $I$. 
A map $w: T\to X$ belongs to $W^{1,2}(T,X)$ if $w\circ\varphi\in W^{1,2}(I,X)$ for
some (and thus any) bilipschitz map $\varphi:I\to T$. Such Sobolev curves
arise for instance as follows.
If $u\in W^{1,2}(D,X)$ and $F:(0,1)^2\to D$ is a bilipschitz embedding, then
for almost every $s\in(0,1)$ the curve $u|_{T_s}$ lies in $W^{1,2}(T_s,X)$ where $T_s=F(\{s\}\times(0,1))$. The definition of $W^{1,2}(T, X)$ naturally extends to the case where $T$ is bilipschitz to $S^1$. Even more generally, if $T\subset\R^2$ is a finite union of bilipschitz curves $T_i$ as above then a map $w: T\to X$ is in $W^{1,2}(T, X)$ if $w$ has a continuous representative and $w|_{T_i}\in W^{1,2}(T_i, X)$ for every $i$.

Now suppose $X$ supports a $(C,l_0)$-isoperimetric inequality
and $J\subset D$ is a bilipschitz Jordan curve enclosing a Jordan domain $\Om$.
Then for every $\ga\in W^{1,2}(J,X)$ of length $l< l_0$
there exists $u\in W^{1,2}(\Om,X)$ with $\a(u)\leq C\cdot l^2$ \cite[Lemma~4.6]{LWint}.

\subsection{Conformal and almost conformal Sobolev discs}
Let $\Om\subset\R^2$ be a domain.
A map $u\in W^{1,2} (\Omega,X)$ is called {\em quasi-conformal} if there exists a constant $Q\geq 1$ such that  at almost all $z\in \Omega$ we have   
\[\apmd u_z(v)\leq Q\cdot \apmd u_z(w)\] for all $v, w\in S^1$. In this case we say that $u$ is {\em $Q$-quasi-conformal}, and if $Q$ can be chosen equal to $1$, we call $u$ {\em conformal}. In  this case, $\apmd u_z$ is a multiple $f(z)\cdot s_0$ of the standard Euclidean norm $s_0$ on $\R^2$. The function
$f\in L^2 (\Omega)$ will be called the {\em conformal factor} of $u$. The conformal factor $f$ of a conformal map $u\in W^{1,2} (\Omega,X)$ coincides with the generalized gradient  $\rho _u$.

For every Sobolev disc $u\in W^{1,2}(D,X)$ 
holds
\[\a_{\mu^i}(u)\leq E_+^2(u).\] 
Moreover, if equality holds then $u$ is $\sqrt{2}$-quasi-conformal \cite[Corollary~3.3]{LW_energy}.

 Recall from \cite{LW_plateau} that a complete metric space $X$ is said to have {\em property (ET)} if for every Sobolev disc in $X$ the approximate metric derivative 
comes from a possibly degenerate inner product at almost every point. If $X$ is a complete metric space with property (ET), then 
 parametrized area and inscribed Riemannian area of a Sobolev disc $u\in W ^{1,2}(D,X)$  
coincide,
\[\a_{\mu^i}(u)=\a(u).\]
Moreover, the equality $\a(u)=E_+^2(u)$ implies that $u$ is conformal.

The following version of Morrey's $\eps$-conformality lemma allows to find good parametrizations of Sobolev discs.

\bt[{\cite[Theorem~1.4]{FW_morrey}}]\label{thm:Morrey-epsilon-conformality-disc}
 Let $X$ be a complete metric space and let $u\in W^{1,2}(D,X)$. Then for every $\eps>0$ there exists a diffeomorphism $\varphi: D\to D$ such that $$E_+^2(u\circ\varphi) \leq \Area_{\mu^i}(u) + \eps.$$
Moreover, there is such a  map $\varphi$ which extends to a diffeomorphism of $\bar{D}$ and is conformal in a neighbourhood of the boundary.
\et
Note that the use of inscribed Riemannian area is essential here.
\subsection{Regularity of quasi-conformal Sobolev discs}\label{sec_reg}

The following interior regularity result is a consequence of \cite[Proposition~8.4]{LW_plateau} and \cite[Proposition~8.7]{LW_plateau}.

We say that a property holds for almost every bilipschitz Jordan curve in $D$,
if whenever $\varphi:S^1\times[0,1]\to D$ is a bilipschitz annulus, then the property holds true for
almost every circle $\varphi(S^1,t)$.

\begin{thm} \label{thm_reg_new}
Let $Z$ be a complete metric space and $u\in W^{1,2}(D,Z)$ a
$Q$-quasi-conformal Sobolev disc.
Suppose that there exist $C, l_0>0$  such that for almost every bilipschitz Jordan curve $\ga$ in $D$ with $\ell(u|_\ga)<l_0$ we have
\begin{equation}\label{eq_balls}
    \a(u|_{\Om_\ga})\leq C\cdot\ell^2(u|_{\ga}),
\end{equation}
where $\Om_\ga$ is the Jordan domain enclosed by $\ga$.
Then the following statements hold:
\begin{enumerate}
    \item There exists $p>2$ such that $u\in W^{1,p}_{loc}(D,Z)$. In particular, $u$ has a continuous representative $\bar{u}$ which moreover satisfies Lusin's property (N). 
    \item The representative $\bar{u}$ is locally H\"older continuous. In fact, for every $\delta\in(0,1)$ there exists $L>0$ such that for all
 $z_1, z_2\in\bar B_\delta(0)$ there exists a path $\ga$
in $\bar B_\delta(0)$ such that 
\[\ell_X(\bar{u}\circ\ga)\leq L\cdot |z_1-z_2|^\alpha,\]
where $\alpha=\frac{1}{4\pi Q^2 C}$.
\item If $\tr u$ has a continuous representative then $\bar{u}$ continuously extends to $\bar{D}$.
\end{enumerate}
\end{thm}

The proof can be assembled from arguments in \cite{LW_plateau}.
It basically follows from \cite[Theorem~8.2]{LW_plateau},
\cite[Proposition~8.7]{LW_plateau} and \cite[Theorem~9.1]{LW_plateau}.
These results assume the map $u$  to be minimal and the space $X$ to
support a $(C,l_0)$-isoperimetric inequality. 
While this differs from our setting, these additional assumptions are only
used in \cite{LW_plateau} to prove \cite[Lemma~8.6]{LW_plateau} which ensures that inequality \eqref{eq_balls} holds.

For the convenience of the reader let us give a more detailed account  on how to obtain Theorem~\ref{thm_reg_new} from the results in \cite{LW_plateau}.
We start by proving (1) which corresponds to  \cite[Proposition~8.4]{LW_plateau}, namely higher integrability and therefore 
H\"older continuity of $u$. Next we prove (2) which corresponds to \cite[Proposition~8.7]{LW_plateau}, the intrinsic H\"older regularity. 
The proofs in \cite{LW_plateau} apply because they only use the quasi-conformality of $u$ and inequality \eqref{eq_balls} for balls entirely contained in $D$.  
We are left with (3) which corresponds to \cite[Theorem~9.1]{LW_plateau}.
To show that $u$ extends continuously to a point $z\in S^1$, the proof in \cite{LW_plateau}
considers a conformal diffeomorphism $\varphi:D\to\Om$ where $\Om=D\cap B_r(z)$
for some $z\in S^1$ and a suitable $r\in(0,1)$. To make the argument from \cite{LW_plateau} work, we need that 
\[\a(u\circ\varphi|_{B_s(x)})\leq C\cdot\ell^2(u\circ\varphi|_{\D B_s(x)})\]
holds for all $x\in D$ and almost all $s<1-|x|$. However, this is guaranteed since $\varphi$
is bilipschitz on every compact subset of $D$ and inequality~\eqref{eq_balls}
holds for almost every bilipschitz Jordan curve by assumption.

\section{Dehn functions, thickenings, and minimal discs}\label{sec_dehn}

\subsection{Filling area and Dehn functions}\label{sec_filling_area}

Recall that for a Jordan curve $\Ga$ in a metric space $X$ we denote by $\La(\Ga,X)$
the family of Sobolev discs $u\in W^{1,2}(D,X)$ whose traces have representatives which are weakly monotone
parametrizations of $\Ga$. We set

\[\Fill_X(\Ga):=\inf\{\a(u)\mid u\in\La(\Ga,X)\}.\]
 Similarly, if $c:S^1\to X$ is a curve, we set  
\[\Fill_X(c):=\inf\{\a(u)\mid u\in W^{1,2}(D,X),\, \tr(u)=c\}.\]
Notice that if $c$ is a weakly monotone parametrization of a Jordan curve $\Ga$,
then, by definition, $\Fill_X(\Ga)\leq\Fill_X(c)$.

The {\em Dehn function} of $X$ is given by
\[\delta_X(r):=\sup\{\Fill_X(c)\mid c:S^1\to X \text{ Lipschitz, } \ell(c)\leq r\}\]
for all $r>0$. Similarly, we define the Riemannian versions 
\[\Fill^{\mu^i}_X(\Ga)\; \text{ and } \;\Fill^{\mu^i}_X(c)\; \text{ and }\delta^{\mu^i}_X(r)\]
by replacing the Hausdorff area of Sobolev discs by the inscribed Riemannian area. Clearly, Dehn functions are non-decreasing in $r$. Notice that $\delta_X(r)\leq \delta_X^{\mu^i}(r)$, and equality holds for example when $X$ has property (ET). 

As mentioned in the introduction, we denote by $\delta_\kappa$ the Dehn function of the model surface $M^2_\kappa$. 
Explicitly, for $r\in(0,2 D_\ka)$ we have 
\[\delta_\ka(r)=
\begin{cases}
    \frac{2\pi}{\ka}-\sqrt{(\frac{2\pi}{\ka})^2-\frac{r^2}{\ka}}    & \quad \text{for  } \ka>0\\
    \frac{r^2}{4\pi}  & \quad \text{for } \ka=0\\
    \sqrt{(\frac{2\pi}{\ka})^2+\frac{r^2}{\ka}}-\frac{2\pi}{\ka} & \quad \text{for } \ka<0
  \end{cases}
\]
as follows for instance from the isoperimetric inequality in $M^2_\ka$, see \cite{O_iso}.
Note that for $\ka>0$ we have the quadratic bound $\delta_\ka(r)\leq\frac{r^2}{2\pi\ka}$ for $r\in(0,2D_\ka)$.

\subsection{Universal thickenings}

Let $X$ and $Y$ be metric spaces and $\epsilon>0$. We say that $Y$ is an \emph{$\epsilon$-thickening} of $X$ if there exists an isometric embedding $\iota: X\to Y$ such that the Hausdorff distance between $\iota(X)$ and $Y$ is at most $\epsilon$.  The embedding $\iota$ is then a $(1,\epsilon)$-quasi-isometry.
 We will make use of the following which can be proved in the same way as \cite[Proposition~3.5]{LWY_dehn}.
 Recall that a metric space $X$ is said to be {\em $L$-Lipschitz 1-connected up to scale $\la_0$}
 for some $L\geq 1$ and $\la_0>0$ if every $\la$-Lipschitz curve $c:S^1\to X$
 with $\la<\la_0$ extends to a $L \la$-Lipschitz map on $\bar D$.
 If the scale $\la_0$ is not important, we simply say $X$ is $L$-Lipschitz 1-connected up
 to some scale. We clearly have $\delta_X(r)\leq\deltalip_X(r)$ for all $r>0$, where $\deltalip_X$ denotes the Lipschitz Dehn function defined by filling Lipschitz curves by Lipschitz discs.  However, equality holds for all complete length spaces which are Lipschitz 1-connected up 
to some scale \cite[Proposition~3.1]{LWY_dehn}. The same applies to the $\mu^i$-versions of the Sobolev and Lipschitz Dehn functions. 

\bprop\label{prop_thickening}
 There exists $L\geq 1$ with the following property. Let $X$ be a complete length space and let $\eps, r_0>0$. Suppose $\delta:(0,r_0)\to\R$ is continuous and $$\delta^{\mu^i}_X(r)\leq \delta(r)+\eps^2$$ for all $r\in(0,r_0)$. Then there exists a complete length space $Y$ which is a $\eps$-thickening of $X$ and $L$-Lipschitz $1$-connected up to scale $\frac{\eps}{L}$ and satisfies $$\delta^{\mu^i}_Y(r)\leq \delta(r)+Lr^2$$ for all $r\in(0,r_0)$ and if $r\in(\sqrt{\eps}, r_0)$ then $$\delta^{\mu^i}_Y(r)\leq \delta(r)+ \sqrt{\eps}r^2.$$
\eprop

\begin{rem}
If $\delta(r)=O(r^2)$ then there exist $C\geq 1$ and $r_1\in(0,r_0)$ depending only on $L$ and the function $\delta$ such that $Y$ has a $(C,r_1)$-isoperimetric inequality. 
\end{rem}

\bdfn
For a length space $X$ we call an $\eps$-thickening $X_\eps$ as in Proposition~\ref{prop_thickening} a {\em universal $\eps$-thickening}.
\edfn

\subsection{Plateau problem and intrinsic minimal discs}\label{sec_int_min}

\bdfn
Let $X$ be a complete metric space and $\Ga\subset X$ a Jordan curve.
We call $u\in\La(\Ga,X)$ a {\em solution to the Plateau problem} for the curve $\Ga$,
if $\a(u)=\Fill(\Ga)$ and $u$ has minimal energy among all area minimizers in $\La(\Ga,X)$.
A solution to the Plateau problem will sometimes simply be called a {\em minimal disc}.
\edfn

In this section we consider the following setting. Let $X$ be a complete length space which satisfies property (ET).
Let $\Ga\subset X$ be a rectifiable Jordan curve and $u\in\La(\Ga,X)$ a solution to the Plateau problem. In particular, $u$ is conformal
\cite[Theorem~11.3]{LW_plateau}. We assume that $u$ satisfies inequality~\eqref{eq_balls} and therefore, by Theorem~\ref{thm_reg_new},  has a (intrinsically) H\"older continuous representative which continuously extends to $\bar D$. Denote this representative still by $u$.
Then \cite[Theorem~1.1]{LWint} yields the following
structure for the intrinsic minimal disc $Z_u$, cf. Section~\ref{sec_associated}. 
The setting in \cite{LWint} asks $X$ to support a $(C,l_0)$-isoperimetric inequality. However, the proof of \cite[Theorem~1.1]{LWint} only uses the quasi-conformality of $u$ and inequality~\eqref{eq_balls}.

\bthm
The intrinsic minimal disc $Z_u$ is a compact geodesic space. The canonical projection $P:\bar D\to Z_u$
is continuous. The map $u:\bar D\to X$ has a canonical factorization $u=\bar u\circ P$,
where $\bar u:Z_u\to X$ is 1-Lipschitz. For any curve $\ga$ in $\bar D$ the lengths of
$P\circ \ga$ and $u\circ\ga$ coincide, thus $\bar u$ preserves the length of $P\circ\ga$.
\ethm

Now suppose $\delta:(0,r_0)\to\R$ is a continuous non-decreasing function and there exists $0< r_1 < r_0$ such that $\ell(\Ga)<r_1$ and for every Jordan domain $\Om\subset D$ with $\ell(u|_{\D\Om})< r_1$ holds
\begin{equation}\label{eq_int_area_u_Om}
    \a(u|_\Om)\leq\delta(\ell(u|_{\D\Om})).
\end{equation}
In this setting, the arguments in \cite{LWint} provide strong topological and isoperimetric properties of the 
intrinsic minimal disc:

\bthm\label{thm_int_top}
The map $P:\bar D\to Z_u$ is a uniform limit of homeomorphisms.
For every Jordan domain $\Om\subset Z_u$  with $\ell(\D\Om)<r_0$ holds
\begin{equation}\label{eq_int_disc}
  \hm^2(\Om)\leq \delta(\ell(\D\Om)).  
\end{equation}
\ethm

The proof can be assembled from \cite{LWint}.
 By \cite[Lemma~6.3, Lemma~6.4, Corollary~4.5, Theorem~8.1]{LWint} the natural projection satisfies the first statement. This part relies on 
Moore's recognition theorem for 2-manifolds \cite[Theorem~7.11]{LWint}.
The second statement follows from \cite[Theorem~8.2]{LWint}, here is where the continuity of $\delta$ is needed.
We emphasize that unlike in inequality \eqref{eq_int_area_u_Om}, Theorem~\ref{thm_int_top} does not require that the length of $\D\Om$ be bounded by $r_1$. Indeed, \cite[Theorem~1.1]{LWint} shows that $\ell(\D Z_u)=\ell(\Ga)<r_1$ and therefore $\hm^2(Z_u)\leq \delta(r_1)$. In particular, even if 
$\Om\subset Z_u$ is a Jordan domain with
$r_1<\ell(\D\Om)<r_0$, we still have 
$\hm^2(\Om)\leq \delta(r_1)\leq\delta(\D\Om)$ by the monotonicity of $\delta$.

\section{Euclidean tangent planes}\label{sec:ET}

Recall once again that a complete metric space $X$ is said to have property (ET) if for every $u\in W^{1,2}(D, X)$ the approximate metric derivative 
$\apmd u_z$ comes from a possibly degenerate inner product at almost every $z\in D$. 
 
 The proof of the following proposition is very similar  to that of \cite[Theorem 3.1]{W_dehn}, which was originally inspired by \cite[Theorem 5.1]{W_gromov}. The statement generalizes \cite[Theorem 5.2]{LWcurv} and \cite[Theorem 3.1]{W_dehn}.

\bprop\label{prop:general-ET}
 Let $r_0>0$ and let $\delta: (0,r_0)\to\R$ be a continuous  non-decreasing function satisfying $$\limsup_{r\to0} \frac{\delta(r)}{r^2}\leq \frac{1}{4\pi}.$$ Let $(X_n)$ be a sequence of complete length spaces satisfying $$\delta_{X_n}(r)\leq (1+\eps_n)\cdot\delta(r) + \eps_n$$ for all $r\in(0, r_0)$, where $\eps_n>0$ tends to zero as $n\to\infty$. Then every ultralimit of $(X_n)$ has property (ET).
\eprop

We provide the proof for the convenience of the reader.
 
\proof
Let $X_\om=(X_\om, d_\om)$ be an ultralimit of the sequence $(X_n)$ and suppose, by contradiction, that $X_\om$ does not have property (ET). By \cite[Lemma~3.2]{W_dehn} there exists a non-Euclidean norm $\|\cdot\|$ on $\R^2$ with the following property. For every finite set $\{v^1,\dots, v^m\}\subset\R^2$ and every $\lambda>1$ there exist points $x^1,\dots, x^m\in X_\om$ and $\eta>0$ arbitrarily small such that 
\begin{equation}\label{eq-property-non-Euclidean-finite}
    \lambda^{-1}\eta\|v^k-v^l\| \leq d_\om(x^k,x^l) \leq \lambda\eta \|v^k-v^l\|
\end{equation}
for all $k, l=1,\dots, m$. We denote by $V$ the normed space $(\R^2, \|\cdot\|)$ and let $\mathbb{I}_V\subset V$ be an isoperimetric subset for $V$, that is, $\mathbb{I}_V$ is convex and has largest area among all convex subsets of $V$ with prescribed boundary length. Since $V$ is non-Euclidean we have  
\begin{equation}\label{eq:non-Eucl-area-bound}
 \hm_V^2(\mathbb{I}_V) > \frac{1}{4\pi}\cdot \ell^2_V(\partial \mathbb{I}_V),
\end{equation}
  see for example \cite[Lemma 5.1]{LWcurv}.
 
  Let $\gamma: S^1\to V$ be a constant speed parametrization of $\D\mathbb{I}_V$. Choose $\lambda>1$ sufficiently close to $1$ and $m\in\N$ sufficiently large, both to be determined later. For $k=1,\dots, m$ set $z_k= e^{2\pi i\frac{k}{m}}\in S^1$ and let $v^k:= \gamma(z_k)$. By the above, there exist $x^1,\dots, x^m\in X_\omega$ and $\eta>0$ arbitrarily small such that \eqref{eq-property-non-Euclidean-finite} holds. If $\eta>0$ is small enough then, after replacing the norm $\|\cdot \|$ by the rescaled norm $\eta \cdot \|\cdot\|$, we may assume that $\eta= 1$, that $r_1:=\lambda^4\ell_V(\partial\mathbb{I}_V)$ satisfies $r_1<r_0$ and $\delta(r)\leq \frac{\lambda}{4\pi}\cdot r^2$ for all $0<r\leq r_1$.

  For $k=1,\dots, m$ write $x^k$ as $x^k = [(x_n^k)]$ with $x_n^k\in X_n$. There exists $N\subset\N$ with $\om(N)=1$ and such that $$\lambda^{-1}d_\om(x^k, x^l)\leq d_n(x_n^k, x_n^l)\leq \lambda\cdot d_\om(x^k, x^l)$$ for all $k, l=1,\dots, m$ and all $n\in N$. Fix $n\in N$ large enough, to be determined later, and let $c: S^1\to X_n$ be a Lipschitz curve such that $c(z_k) = x_n^k$ and such that $c$ is a $(\lambda-1)$-geodesic on the segment of $S^1$ between $z_k$ and $z_{k+1}$ for each $k$. It follows from the above that 
$$\ell(c)\leq \lambda^3\cdot \ell_V(\partial\mathbb{I}_V).$$
 
\begin{equation*}
        \a(u)\leq (1+\eps_n)\cdot \delta(\ell(c)) + 2\eps_n\leq \frac{1+\eps_n}{4\pi}\cdot \lambda^7\ell^2_V(\partial\mathbb{I}_V) +2\eps_n.
\end{equation*}
View $V$ as a linear subspace of the space $\ell^\infty$ of bounded sequences, equipped with the supremum norm. Since $\ell^\infty$ is an injective metric space there exists a $\lambda^2$-Lipschitz map $\varphi: X_n\to \ell^\infty$ which maps $x_n^k$ to $v^k$ for all $k$. It follows that the map $\varphi\circ u$ is Sobolev with $\a(\varphi\circ u)\leq \lambda^4\a(u)$. Hence, by \cite[Proposition 3.1]{LWY_dehn}, there exists a Lipschitz map $v:\bar D\to \ell^\infty$ with $v|_{S^1}= \varphi\circ c$ and 
$$\a(v) \leq \a(\varphi\circ u)+ \eps_n\leq \frac{1+\eps_n}{4\pi}\cdot \lambda^{11}\cdot\ell^2_V(\partial\mathbb{I}_V)+2\lambda^4\eps_n +\eps_n.$$

Finally, one constructs exactly as in the proof of \cite[Theorem~3.1]{W_dehn} a Lipschitz homotopy $\varrho: S^1\times[0,1]\to\ell^\infty$ between $\varphi\circ c$ and $\gamma$ with $$\a(\varrho)\leq \frac{C(1+\lambda^3)^2}{m}\cdot\ell^2_V(\partial\mathbb{I}_V),$$ where $C$ is a universal constant. Gluing $\varrho$ and $v$ we obtain a Lipschitz map $w: \bar D\to \ell^\infty$ with $w|_{S^1}=\gamma$ and such that 
\begin{equation}\label{eq:bound-area-w}
 \a(w)\leq \left(\frac{1+\eps_n}{4\pi}\cdot \lambda^{11}+\frac{C(1+\lambda^3)^2}{m}\right) \cdot\ell^2_V(\partial\mathbb{I}_V) +2\lambda^4\eps_n +\eps_n.   
\end{equation}
 Since $\hm^2_V(\mathbb{I}_V)\leq \a(w)$ by the quasi-convexity of the Hausdorff $2$-measure, see \cite{BI_min}, inequality \eqref{eq:bound-area-w} clearly contradicts inequality \eqref{eq:non-Eucl-area-bound} for $\lambda>1$ sufficiently close to $1$, and $m$ and $n$ sufficiently large. This completes the proof.
\qed

\section{Filling area in ultralimits}\label{sec_fill}

Let $C, r_1>0$ and for $n\in\N$ let $X_n$ be a complete length space which is $C$-Lipschitz $1$-connected up to some scale and admits a $(C, r_1)$-isoperimetric inequality. Suppose furthermore that 
$$\delta^{\mu^i}_{X_n}(r)\leq (1+\eps_n)\cdot\delta(r)+ \eps_n r^2$$ for all $r\in(\eps_n, r_0)$, where $\delta:(0,r_0)\to\R$ is continuous and non-decreasing with $$\limsup_{r\to0}\frac{\delta(r)}{r^2} \leq \frac{1}{4\pi}$$ and $\eps_n\to 0$. 
Let $X_\omega$ be the ultralimit of $(X_n)$ with respect to some basepoints and a non-principal ultrafilter $\om$. The following proposition plays a key role in the proofs of our main theorems.

\bprop\label{prop_ultrafillarea}
Let $\Ga_\om\subset X_\om$ be a rectifiable Jordan curve which is the ultralimit of a sequence of rectifiable Jordan curves
$\Ga_n\subset X_n$ of uniformly bounded length.
Then  
\[\wlim \Fill^{\mu^i}_{X_n}(\Ga_n)\leq\Fill_{X_\om}(\Ga_\om).\]
\eprop

Recall that we always have $\Fill_{X_n}(\Ga_n)\leq \Fill^{\mu^i}_{X_n}(\Ga_n)$ and thus 
$$\wlim \Fill_{X_n}(\Ga_n)\leq\Fill_{X_\om}(\Ga_\om)$$ holds as well.
A similar remark applies to the next proposition.

Before we turn to the proof, we begin with a version for parametrized curves.

\bprop\label{prop_lift_fillingarea}
 Let $(c_n)$ be a bounded sequence of Lipschitz curves $c_n: S^1\to X_n$ with uniformly bounded Lipschitz constants. Let $c=\wlim c_n$ be the ultralimit of this sequence. Then \[\wlim\Fill^{\mu^i}_{X_n}(c_n) \leq \Fill_{X_\omega}(c).\]
\eprop

This is a variant of \cite[Theorem~5.1]{W_dehn} and the proof therein applies with minor modifications.

\begin{figure}
    \centering
\includegraphics[scale=0.2,trim={-3cm -2cm 0cm -2cm},clip]{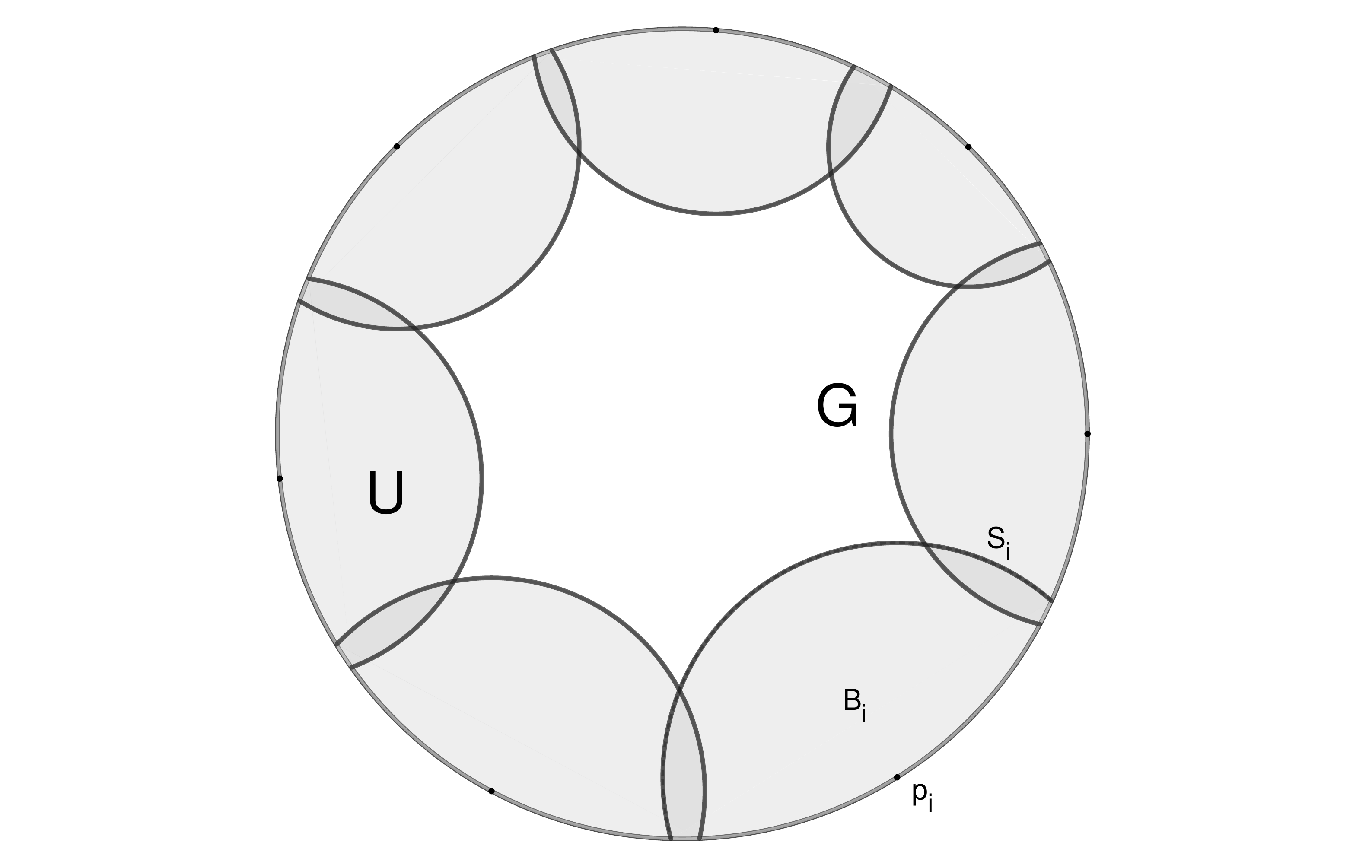}
    \caption{A framed collar $U$ with frame $G$.}
    \label{fig:framed-collar}
\end{figure}

Note that there is no obvious way how to deduce Proposition~\ref{prop_ultrafillarea} from Proposition~\ref{prop_lift_fillingarea}. In the presence of a (local) quadratic isoperimetric inequality one can relate the filling area of any parametrization of a 
Jordan curve to the filling area of its arclength parametrization \cite[Lemma~4.8]{LWint}. 
However,  we do not know whether $X_\om$ admits a (local) quadratic isoperimetric inequality. In particular, the filling area of any Lipschitz parametrization of $\Ga_\om$ might be much smaller than the filling area of a general parametrization.
To overcome this difficulty we need some preparation. In particular, the following notion of framed collar will be useful; see Figure~\ref{fig:framed-collar} for an illustration.

\bdfn
A {\em framed collar} $U\subset\bar D$ is a finite union of closed balls $B_i=\bar B_{r_i}(p_i)\cap \bar{D}$, $i=1,\ldots, m$, centered at points $p_i\in S^1$ such that the open balls cover $S^1$ and non-consecutive balls are disjoint. In particular, $U$ is a  topological annulus whose boundary circles are given by $\D^- U=S^1$
and a bilipschitz Jordan curve $\D^+ U$.
The {\em frame} $G$ of $U$ is the finite graph given by   
$$G=S^1\cup\bigcup\limits_{i=1}^m S_i,$$ where $S_i=\D B_i$. The set $U\setminus G$
is a disjoint union of {\em complementary (open) discs} $\Om_j$, $j=1,\ldots, 2m$.
\edfn

Let us make two simple comments which will be used implicitly later on.
First, if a collar $U$ is contained in the $\eps$-tubular neighborhood of $S^1$ for some $\eps>0$, then every complementary disc $\Om_j$ has diameter at most $2\eps$.
Secondly,  if $U$ has $2m$ complementary discs, then for every small enough $t\in(0,1)$ the intersection $\D B_{1-t}(0)\cap G$ consists of $2m$ points.

We will use framed collars in combination with the isoperimetric inequality to produce Sobolev homotopies of small area.

\blem\label{lem_graph}
Let $X$ be a complete metric space and $u\in W^{1,2}(D,X)$.
Then for every $\eps>0$ there  exists a framed collar
$U\subset N_\eps(S^1)$ in $\bar D$ with complementary discs $\Om_i$ such that   
\begin{itemize}
 \item $u|_{G\cap D}\in W^{1,2}(G\cap D,X)$ and
	\[\sum\limits_{i=1}^{2m} \ell^2(u|_{\D\Om_i\cap D})<\eps;\] 
\item  $u|_{\D^+ U}=\tr(u|_{D\setminus U})$.
\end{itemize}
Moreover, $G$ can be chosen to omit a given $2$-exceptional family of curves in $D$.
\elem

The proof relies on a construction used in the proof of \cite[Lemma~4.8]{LWint}.

\proof
We may assume that $u$ is absolutely continuous on $2$-a.e.~curve in $D$, see Section~\ref{sec_length}.
Fix $\eps>0$ and let $\delta>0$ be a constant whose size will be determined in terms of $\eps$.
We choose a small $\rho= \sin (\frac {2\pi} m)<\eps$ such that the restriction  of $u$ to the $\rho$-neighborhood of $S^1$ in 
$D$  has energy at most $\delta$. 
Next, we choose  equidistant points $p_1,....,p_m$ on $S^1$ with pairwise Euclidean distance $\rho<\eps$.

Denote by $E_i$ the energy of the restriction of $u$ to $B_\rho(p_i)\cap D$. Note that 
\[\sum\limits_{i=1}^m E_i\leq 2\cdot E^2_+(u|_{N_\rho(S^1)})\leq2\delta.\]
By \cite[Lemma~3.5]{LWint}, we find   subsets  $R_i\subset (\frac 2 3 \rho , \rho)$ of
positive measure, such that the following holds true for every $r_i\in R_i$.   
The restriction of   $u$  to  the  distance circle $S_i$  of radius $r_i$ around
$p_i$ in $D$ is a continuous curve in $W^{1,2} (S_i,X)$ and its length $\ell_i$ 
satisfies $\ell_i^2<6\pi \cdot  E_i$.
We define $B_i$ to be the ball $\bar B_{r_i}(p_i)$. By construction,   
$U=\bigcup\limits_{i=1}^m\bar B_i$ is a framed collar with frame
$G:=S^1\cup\bigcup\limits_{i=1}^m S_i$.
The domain $U$ is subdivided by the circular arcs  $S_i$ in $2m$ Lipschitz discs $\Om_j$. The boundary of any $\Om_j$
 consists of two or three
parts  of consecutive circles $S_i$ and a part $\D\Om_j\cap S^1$. 
By construction, $u|_{G\cap D}\in W^{1,2}(G\cap D,X)$ and
\[\sum\limits_{j=1}^{2m} \ell^2(u|_{\D\Om_j\cap D})\leq 6\cdot\sum\limits_{i=1}^{m}\ell^2_i\leq 36\pi\cdot\sum\limits_{i=1}^{m}E_i\leq 72\pi\cdot\delta.\]
Hence we can choose $\delta=\delta(\eps)$ small enough to guarantee the claimed length bound.

The last two statements hold since we can choose $r_i$ freely from the positive measure set $R_i$.
\qed

\proof[Proof of Proposition~\ref{prop_ultrafillarea}]
Fix $\eps_0>0$.
We may assume that  $\La(\Ga_\om,X_\om)$ is non-empty, otherwise there is nothing to show.
Choose $u\in\La(\Ga_\om,X_\om)$ such that $\a(u)<\Fill_{X_\om}(\Ga_\om)+\frac{\eps_0}{4}$.
We choose $\eps>0$ small enough such that $\a(u|_{N_\eps(S^1)})<\frac{\eps_0}{4}$. 
By Lemma~\ref{lem_graph}, we find a framed collar 
$U\subset N_\eps(S^1)$ with frame $G$ and
complementary discs $\Om_i\subset U$ such that 
$u|_{G\cap D}\in W^{1,2}(G\cap D,X)$ with
	\[\sum\limits_{i=1}^l \ell^2(u|_{\D\Om_i\cap D})<\eps.\] 
Since $G$ can be chosen to omit a given $2$-exceptional family of curves in $D$, we may furthermore assume that the curves $u|_{\D\Om_i\cap D}$ and $\tr(u)|_{\D\Om_i\cap S^1}$ form a closed loop.
By adjusting $\eps$ if necessary, we can therefore arrange
\[\sum\limits_{i=1}^l \ell^2(f|_{\D\Om_i})<\frac{\eps_0}{16C},\]
where $f\in W^{1,2}(G,X_\om)$ denotes the map given by $u$ on $G\cap D$ and a constant speed parametrization of $\tr(u)|_{\D\Om_i\cap S^1}$ on every interval $\D\Om_i\cap S^1$.
Let $\bar f:G\to X_\om$ be the Lipschitz reparametrization of $f$ which has constant speed on every edge of $G$.
By \cite[Lemma~2.6]{LWY_dehn} we can lift $\bar f$  to Lipschitz graphs $\bar f_{n}:G\to X_n$ with uniformly controlled Lipschitz constant such that  
$\ell(\bar f_{n})\leq 2\cdot\ell(\bar f)$.
Since producing such lifts only involves extensions from finite sets, we can arrange that for every $n\in\N$
the $\bar f_{n}$ restrict on $S^1$ to a Lipschitz parametrization of $\Ga$. Set $\eta=\D^+ U$. Then $\eta$ is a bilipschitz Jordan curve in $D$
and $(\bar f_{n}|_\eta)$ is a bounded sequence of Lipschitz curves with uniformly controlled Lipschitz constants and, by construction, 
\[\bar f|_\eta=\wlim \bar f_{n}|_\eta.\]
By Proposition~\ref{prop:general-ET}, $X_\om$ satisfies property (ET), thus \cite[Lemma~2.6]{LW_harm} and
Proposition~\ref{prop_lift_fillingarea} imply that there exists $N\subset\N$ with $\om(N)=1$ such that  
 \[\Fill^{\mu^i}_{X_n}(\bar f_n|_\eta)\leq\Fill_{X_\om}(u|_\eta)+\frac{\eps_0}{4}\]
for every $n\in N$.
By our choice of $u$ and $\eps$, we have
\[\Fill_{X_\om}(u|_\eta)\leq\a(u)+\frac{\eps_0}{4}\leq\Fill_{X_\om}(\Ga_\om)+\frac{\eps_0}{2}.\]
On the other hand, we can fill the curves $\bar f_n|_{\D\Om_i}$ using the $(C,r_1)$-isoperimetric inequality to produce
Sobolev annuli $h_n\in W^{1,2}(A,X_n)$ which join $\bar f_n|_\eta$ and $\bar f_n|_{S^1}$ with
\[\a_{\mu^i}(h_n)\leq C\cdot\sum_{i=1}^l \ell^2(\bar f_n|_{\D\Om_i})< \frac{\eps_0}{4}.\]
Thus, since $\bar f_n|_{S^1}$ is a Lipschitz parametrization of $\Ga\subset X_n$ and $X_n$
supports an isoperimetric inequality, \cite[Lemma~4.8]{LWint} implies 
\[|\Fill^{\mu^i}_{X_n}(\bar f_n|_\eta)-\Fill^{\mu^i}_{X_n}(\Ga)|<\frac{\eps_0}{4}.\]
Putting things together, we obtain
\[\Fill^{\mu^i}_{X_n}(\Ga)\leq\Fill_{X_\om}(\Ga_\om)+\eps_0\]
for all large enough $n\in N$.
\qed

\section{Intrinsic isoperimetric inequality}\label{sec_int}

Let $r_0, r_1>0$, $C\geq 1$ and $\beta_n\to 0$. Let
$\delta:(0,r_0)\to \R$ be a continuous non-decreasing  function.
Suppose that $(Y_n)$ is a sequence of complete length spaces, each admitting a 
$(C,r_1)$-isoperimetric inequality, and such that  
$$\delta_{Y_n}(r)\leq (1+\beta_n)\cdot\delta(r)+ \beta_n r^2$$
if $r\in(\beta_n, r_0)$.

\bprop\label{prop_inisin}
Suppose that all $Y_n$ are isometrically embedded in a single  complete metric space $Z$.
Let $v\in W^{1,2}(D,Z)$ be a continuous map which admits a continuous extension $v:\bar D\to Z$. 
Further, let $(v_n)$ be a sequence in $W^{1,2}(D,Y_n)$ such that the following holds true:
\begin{itemize}
	\item $v_n\to v$ pointwise almost everywhere;
	\item $\a_{\mu^i}(v_n|_\Om)\to \a_{\mu^i}(v|_\Om)$ for every Jordan domain $\Om\subset D$ with $v|_{\D\Om}$ rectifiable;
	\item $\ell(v_n|_\gamma)\to \ell(v|_\gamma)$ for 2-a.e.~curve $\ga$ in $D$;
	\item $\a_{\mu^i}(v_n)\leq\Fill^{\mu^i}_{Y_n}(\tr(v_n))+\eps_n$ for some $\eps_n\to 0$.
\end{itemize}
Let $\Om\subset D$ be a Jordan domain such that $v|_{\D\Om}$ is rectifiable
with $\ell(v|_{\D\Om})<r_0$.
Then 
\[\a_{\mu^i}(v|_{\Om})\leq \delta(\ell(v|_{\D\Om})).\]
\eprop

\proof
Let $F:D\to\Om$ be a conformal diffeomorphism and set $u_n:=v_n\circ F$ and $u:=v\circ F$. 
Recall that by Caratheodory's theorem, $F$ extends to a homeomorphism $\bar D\to\bar\Om$.
Since $v$ is continuous on $\bar D$, the map $u:\bar D\to Z$ is continuous and  lies in $W^{1,2}(D,Z)$.

Since $F$ preserves $2$-exceptional families of curves \cite[Proposition~3.5]{LWcan}, we have 
$\ell(u_n|_\gamma)\to \ell(u|_\gamma)$ for 2-a.e.~curve $\ga$ in $D$.

 By Lemma~\ref{lem_graph},
for every $\eps>0$ there  exists a 
framed collar $U\subset N_\eps(S^1)$ with frame $G$. Moreover, $u|_{G\cap D}\in W^{1,2}(G\cap D,Z)$ and
	\[\sum\limits_{i=1}^l \ell^2(u|_{\D\Om_i\cap D})<\eps\]
	holds. In addition, $u|_{\D^+ U}=\tr(u|_{D\setminus U})$.
Since $G$ can be chosen to omit a given $2$-exceptional family of curves in $D$, we may additionally assume:
\begin{itemize}
	\item $u_n|_{G\cap D}\in W^{1,2}(G\cap D,Z)$ for all $n\in\N$;
	\item $u_n|_{G\cap D}\to u|_{G\cap D}$ pointwise $\hm^1$-almost everywhere;
	\item $\ell(u_n|_{G\cap D})\to \ell(u|_{G\cap D})$;
\end{itemize}
By the area convergence assumption and the semi-continuity of area (Proposition~\ref{prop_semi}), we may also assume $\a_{\mu^i}(u_n|_U)<\eps$
for almost all $n\in\N$.

To avoid issues when gluing Sobolev maps we will approximate $\Om$ from within by Lipschitz domains
   as follows.
For $\rho\in(0,1)$ consider the sphere $\si_\rho:=\D B_{1-\rho}(0)\subset D$. For almost all small $\rho\in(0,1)$, the following holds:
\begin{itemize}
	\item The set  
$\si_\rho\cap G$ consits of $l$ points;
\item we have pointwise convergence 
$u_n\to u$ on $\si_\rho\cap G$.
\end{itemize}

Note that $F(\si_\rho)\subset\Om$ is a bilipschitz Jordan curve for every $\rho>0$.

Set $\tilde\Om_i:=\Om_i\cap B_{1-\rho}(0)$.
Now let $f_n:G\to Y_n$ be any continuous map  which is an $\eps$-geodesic on every topological interval $\Om_i\cap \si_\rho$ and which is given by $u_n$ 
on $G\cap B_{1-\rho}(0)$. Then  $f_n\in W^{1,2}(G,Y_n)$ and since the diameter of $\Om_i$ can be assumed to be arbitrary small, 
after possibly adjusting $\eps$ and $U$, we may assume
\[\sum\limits_{i=1}^l \ell^2(f_n|_{\D\tilde\Om_i})<\eps\]
for almost all $n\in\N$. 

Using the $(C,r_1)$-isoperimetric inequality to fill the curves $f_n|_{\D\tilde\Om_i}$, we obtain homotopies $h_n\in W^{1,2}(A, Y_n)$ between $v_n|_{\D^+ U}$ and $f_n|_{\si_\rho}$ with
\[\a_{\mu^i}(h_n)\leq C\cdot\eps. \]

Note that $h_n|_{\D^-A}$ is a piecewise $\eps$-geodesic in $Y_n$ whose vertices are the images under $u_n$
of the points $\{\theta_1,\ldots\theta_l\}:=\si_\rho\cap G$. 
By the pointwise convergence $u_n\to u$ on $\si_\rho\cap G$, we have
$d_{Y_n}(u_n(\theta_i),u_n(\theta_{i+1}))\to d_Z(u(\theta_i),u(\theta_{i+1}))$.
Thus for $n\in\N$ large enough, we have 
\begin{align*}
\ell(h_n|_{\D^-A})&\leq(1+\eps)\sum\limits_{i=1}^{l-1}d_{Y_n}(u_n(\theta_i),u_n(\theta_{i+1}))\\
&\leq(1+\eps)^2\sum\limits_{i=1}^{l-1}d_{Z}(u(\theta_i),u(\theta_{i+1}))\leq(1+\eps)^2\ell(v|_{\D\Om}).
\end{align*}
On the other hand, for small enough $\eps>0$ and large enough $n\in\N$ we have
\[\ell(h_n|_{\D^-A})\geq(1-\eps)\sum\limits_{i=1}^{l-1}d_{Z}(u(\theta_i),u(\theta_{i+1}))\geq\beta_n.\]
Thus, by assumption, we find a filling $w_n\in W^{1,2}(D,Y_n)$ of $h_n|_{\D^-A}\in W^{1,2}(S^1,Y_n)$ with 
\begin{align*}
\a_{\mu^i}(w_n)&\leq (1+\beta_n)\cdot\delta(\ell(h_n|_{\D^-A}))+\beta_n\ell^2(h_n|_{\D^-A})\\
&\leq(1+\beta_n)\cdot\delta((1+\eps)^2\ell(v|_{\D\Om}))+(1+\eps)^4\beta_n\ell^2(v|_{\D\Om}),
\end{align*}

We glue the maps $v_n|_{(D\setminus\Om)\cup F(U)}$, $h_n$ and $w_n$ to obtain a new map $\tilde v_n\in W^{1,2}(D,X_n)$ with 
$\tilde v_n|_{(D\setminus\Om)\cup F(U)}=v_n|_{(D\setminus\Om)\cup F(U)}$. Here we use that $\D^+ U$ and $F(\si_\rho)$ are bilipschitz Jordan curves, 
cf. \cite[Lemma~4.6]{LWint}.
By the last item of our assumptions and our construction so far, we obtain for large enough  $n\in\N$ the estimates
\begin{align*}
\a_{\mu^i}(v_n)&\leq \a_{\mu^i}(\tilde v_n)+\eps_n\\
&\leq\a_{\mu^i}(v_n|_{D\setminus\Om})+\a_{\mu^i}(v_n|_{F(U)})\\
&\quad\quad+\a_{\mu^i}(h_n)+\a_{\mu^i}(w_n)+\eps_n\\
&\leq\a_{\mu^i}(v_n|_{D\setminus\Om})+ (1+\eps)\cdot\delta((1+\eps)^2\ell(v|_{\D\Om}))+2C\eps.
\end{align*}
Thus,
\[\a_{\mu^i}(v_n|_\Om)\leq (1+\eps)\cdot\delta((1+\eps)^2\ell(v|_{\D\Om}))+2C\eps,\]
and since $\a_{\mu^i}(v_n|_\Om)\to\a_{\mu^i}(v|_\Om)$ we obtain the claim.
\qed
\medskip

\section{Solutions to Plateau's Problem in ultralimits}\label{sec_plateau}

Let $C, r_1>0$ and for $n\in\N$ let $Y_n$ be a complete length space which is $C$-Lipschitz $1$-connected up to some scale and admits a $(C, r_1)$-isoperimetric inequality. Suppose furthermore that 
$$\delta^{\mu^i}_{Y_n}(r)\leq (1+\beta_n)\cdot\delta(r)+ \beta_n r^2$$ for all $r\in(\beta_n, r)$, where $\delta:(0,r_0)\to\R$ is continuous and non-decreasing with $$\limsup_{r\to0}\frac{\delta(r)}{r^2} \leq \frac{1}{4\pi}$$ and $\beta_n\to 0$.  

Let $Y_\om$ be the ultralimit of $(Y_n)$ with respect to some basepoints and a non-principal ultrafilter $\om$, and let $\Ga_\om\subset Y_\om$
be a rectifiable Jordan curve of length strictly less than $r_0$.
We then have the following generalization of Theorem~\ref{thm_minimal}.

\bthm\label{thm_plateau}
There is a subsequence of $(Y_n)$ whose ultralimit $\hat Y_\om$ (with respect to the same basepoints and ultrafilter) contains an isometric copy of $\Ga_\om$ and has the following property.
There exists a continuous map $v:\bar D\to \hat Y_\om$
which is a solution of the Plateau problem for $\Ga_\om$ in $\hat Y_\om$ and satisfies
\[\a(v|_\Om)\leq \delta(\ell(v|_{\D\Om}))\]
for every Jordan domain $\Om\subset D$ such that  
$\ell(v|_{\D\Om})<r_0$.
\ethm

The proof will occupy the rest of the section. 
By Lemma~\ref{lem_Jordan_lift} there exist rectifiable Jordan curves $\Ga_n\subset Y_n$ such that $\wlim \Ga_n = \Ga_\om$ and $\ell(\Ga_n) \leq 2\ell(\Ga_\om)$ for all $n\in\N$. Let $u_n:D\to Y_n$ be Sobolev discs $u_n\in\La(\Ga_n,Y_n)$ with $\a_{\mu^i}(u_n)\leq\Fill^{\mu^i}_{Y_n}(\Ga_n)+\frac{1}{2n}$.
By Morrey's $\eps$-conformality lemma \cite[Theorem~1.4]{FW_morrey}, we may assume that 
\[E_+^2(u_n)\leq\a_{\mu^i}(u_n)+\frac{1}{2n}\]
for every $n\in\N$.

After passing to a subsequence, we may assume that the $\Ga_n$ converge in the Gromov-Hausdorff sense to $\Ga_\om$.
 By the compactness theorem \cite[Theorem~3.1]{GW_noncompact},
 after possibly passing to a subsequence, there exist a complete metric space $Z$, isometric embeddings $\varphi_n: X_n \hookrightarrow Z$, a compact subset $K\subset Z$ and $v\in W^{1,2}(D, Z)$ such that $\varphi_n(\Ga_n)\subset K$  for all $n\in \N$, and $v_n:=\varphi_n\circ u_n$ converges to $v$ in $L^2(D, Z)$. 
After passing to a further subsequence, we can ensure that $\varphi_n(\Ga_n)$ Hausdorff converges to a set $\Ga_Z\subset K$.
Note that $\Ga_Z$ has to be isometric to $\Ga_\om$.
By \cite[Lemma~6.7]{W_dehn}, after precomposing with conformal diffeomorphisms of the disc, we may further assume that
$(\tr(v_n))$ uniformly converges to a weakly monotone parametrization of $\Ga_Z$.

After passing to a further subsequence, we may assume $v_n\to v$ pointwise on a full measure subset $M\subset D$ as well as 
$\rho_{v_n}\to\rho$ weakly in $L^2(D)$ where $\rho_{v_n}\in L^2(D)$ is the minimal weak upper gradient of $v_n$ and
$\rho\in L^2(D)$ is a weak upper gradient of $v$ (Lemma~\ref{lem_lim}).

We can now embed $v(M)$ isometrically into the ultralimit $\hat Y_\om:=\wlim Y_n$ and
therefore $v\in W^{1,2}(D,\hat Y_\om)$. Note that possibly  $\hat Y_\om\neq Y_\om$ since we passed to subsequences
several times.
 Recall that since $\hat Y_\om$ has property (ET) by Proposition~\ref{prop:general-ET}, 
the Hausdorff area and the Riemannian inscribed area coincide, $\a(v)=\a_{\mu^i}(v)$.
Our next goal is:
\bprop\label{prop_convergence}
In the setting above, there exists a subsequence $(v_{n_l})$ such that we have
\[\lim\limits_{l\to\infty} \a_{\mu^i}(v_{n_l})= \a(v), \quad \lim\limits_{l\to\infty} E_+ ^2(v_{n_l})= E_+ ^2(v),\]
and 
\[E_+ ^2(v)=\a(v).\]
 Moreover,
\[\a(v)=\Fill_{X_\om}(\Ga_\om)=\lim\limits_{l\to\infty}\Fill^{\mu^i}_{Y_{n_l}}(\Ga_{n_l}).\]
In particular, $v$ is conformal and a solution to the Plateau problem for $\Ga_\om$ in $\hat Y_\om$.
In addition,  
$\rho_{v_{n_l}}\to\rho$ in $L^2(D)$ and $\rho$ is the minimal weak upper gradient for $v$, that is, $\rho=\rho_v$.
\eprop

\proof
Passing to a subsequence, we may assume that $(\Fill^{\mu^i}_{Y_n}(\Ga_n))$ converges.
From Proposition~\ref{prop_ultrafillarea},  we obtain 
\[\lim\limits_{n\to\infty}\Fill^{\mu^i}_{Y_n}(\Ga_n)\leq\Fill_{\hat Y_\om}(\Ga_\om). \]
By construction, we have $\a_{\mu^i}(v_n)=\a_{\mu^i}(u_n)$ and 
\[\Fill^{\mu^i}_{Y_n}(\Ga_n)\leq\a_{\mu^i}(v_n)\leq E_+ ^2(v_n)\leq\a_{\mu^i}(v_n)+\frac{1}{2n}\leq\Fill^{\mu^i}_{Y_n}(\Ga_n)+\frac{1}{n}.\]
In particular, we have 
\[\lim\limits_{n\to\infty}\Fill^{\mu^i}_{Y_n}(\Ga_n)=\lim\limits_{n\to\infty}\a_{\mu^i}(v_n)=\lim\limits_{n\to\infty} E_+ ^2(v_n).\]
Thus, by semi-continuity of area we obtain
\[\Fill_{\hat Y_\om}(\Ga_\om)\leq\a(v)\leq\lim\limits_{n\to\infty}\a_{\mu^i}(v_n)=\lim\limits_{n\to\infty}\Fill^{\mu^i}_{Y_n}(\Ga_n)\leq\Fill_{\hat Y_\om}(\Ga_\om).\]
Therefore, equality holds throughout.
Semi-continuity of energy then yields
\[E_+ ^2(v)\leq\lim\limits_{n\to\infty} E_+ ^2(v_n)=\lim\limits_{n\to\infty}\a_{\mu^i}(v_n)=\a(v)\leq E_+ ^2(v).\]
Thus,
\[\lim\limits_{n\to\infty} E_+ ^2(v_n)=E_+ ^2(v)=\a(v).\]
Since $\rho_{v_n}\to\rho$ weakly in $L^2(D)$ and $E^2_+(v_n)=\int_D \rho_{v_n}^2(z)\, dz$,
we obtain $\rho_{v_{n}}\to\rho$ strongly in $L^2(D)$ and therefore $\rho$ is the minimal weak upper gradient of $v$. Finally, the conformality of $v$ follows from the equality $E^2_+(v)=\a(v)$
and the fact that $\hat Y_\om$ has property (ET).
\qed

\medskip

We will now continue to investigate the properties of $v$. To make use of the convergence $v_n\to v$ we will view $v$ as a map with values in $Z$.
Let us assume that we have already passed to a subsequence as provided by Proposition~\ref{prop_convergence}. In particular, we have the $L^2$-convergence $\rho_{v_n}\to\rho_v$.
This allows us to apply
 Fuglede's lemma \cite{HKST15}. Thus, after passing to a further subsequence, we have
\[\lim\limits_{n\to\infty}\int_\ga|\rho_{v_n}-\rho_v|\ ds=0\]
for $2$-a.e.~curve $\ga$ in $D$. 
Since $v$ is conformal by Proposition~\ref{prop_convergence}, we have 
\[\ell(v\circ\ga)=\int_\ga\rho_v\ ds\]
for $2$-a.e.~curve $\ga$ in $D$.
It follows that the lengths of almost all curves converge:
\[\limsup \ell(v_n\circ\ga)\leq\lim\limits_{n\to\infty}\int_\ga \rho_{v_n}\ ds=\int_\ga \rho_{v}\ ds=\ell(v\circ\ga)\leq\liminf \ell(v_n\circ\ga).\]
The last inequality holds since the well-known lower semi-continuity of length with respect to pointwise convergence of curves 
extends to the setting of almost everywhere pointwise convergence.

Recall that every space $Y_n$ admits a $(C,r_1)$-isoperimetric inequality.

\blem
In the above setting the map
$v$ satisfies inequality~\eqref{eq_balls} with constant $C$:
for almost every bilipschitz Jordan curve $\ga$ in $D$ with $\ell(u|_\ga)<r_1$ holds
\begin{equation*}
    \a(u|_{\Om_\ga})\leq C\cdot\ell^2(u|_{\ga})
\end{equation*}
where $\Om_\ga$ is the Jordan domain enclosed by $\ga$.
\elem

\proof
Suppose that $\ga\subset D$ is a bilipschitz Jordan curve with Jordan domain $\Om_\ga$ and $\ell(v|_\ga)<r_1$ and such that 
$\ell(v_n|_\ga)\to\ell(v|_\ga)$. Our choice of $v_n$ ensures
\[\a_{\mu^i}(v_n|_{\Om_\ga})\leq C\cdot\ell^2(v_n|_\ga)+\frac{1}{2n}\]
for all $n$ large enough. Semi-continuity of area, Proposition~\ref{prop_semi}, yields
\[\a(v|_{\Om_\ga})=\a_{\mu^i}(v|_{\Om_\ga})\leq C\cdot\ell^2(v|_\ga).\]
\qed

 Now Theorem~\ref{thm_reg_new}
 implies that the Sobolev disc $v$  has a locally H\"older continuous representative 
 which continuously extends to $\bar D$ and which we will still denote by $v$.

We next claim that for every Jordan domain $\Om\subset D$ such that $v|_{\D\Om}$ is rectifiable we have 
\begin{equation}\label{eq:area-conv}
    \a_{\mu^i}(v_n|_\Om)\to \a_{\mu^i}(v|_\Om).
\end{equation}
Indeed, since $v|_{\D\Om}$ is rectifiable and hence $\Ha^2(v(\D\Om))=0$ it follows from \cite[Proposition 4.3]{LW_plateau} that $\a_{\mu^i}(v|_{\D\Om})=0$.
Hence, by semi-continuity of area, Proposition~\ref{prop_semi}, any subsequence $v_{n_l}$ satisfies
\begin{align*}
&\a_{\mu^i}(v)=\a_{\mu^i}(v|_\Om)+\a_{\mu^i}(v|_{D\setminus\bar\Om})\\
&\leq\liminf \a_{\mu^i}(v_{n_l}|_\Om)+\liminf\a_{\mu^i}(v_{n_l}|_{D\setminus\bar\Om})\\
&\leq\liminf \left(\a_{\mu^i}(v_{n_l}|_\Om)+\a_{\mu^i}(v_{n_l}|_{\D\Om})+\a_{\mu^i}(v_{n_l}|_{D\setminus\bar\Om})\right)\\
&=\a_{\mu^i}(v),
\end{align*}
which proves \eqref{eq:area-conv}.

Now we can apply Proposition~\ref{prop_inisin} to conclude that  for every  Jordan domain $\Om\subset D$ such that $v|_{\D\Om}$ has length strictly 
less then $r_0$ we have
\[\a(v|_\Om)=\a_{\mu^i}(v|_\Om)\leq\delta(\ell(v|_{\D\Om})).\]

This complete the proof of Theorem~\ref{thm_plateau}.

\section{Main Applications}\label{sec_main}

In this final section we provide the proofs of our main results.
We start with a generalization of Theorem~\ref{thm_intro_stable}:

\bthm\label{thm_stable}
Let $\kappa\in\R$ and let $(X_n)$ be a sequence of complete length spaces
such that the Riemannian Dehn function satisfies 
\[\delta^{\mu^i}_{X_n}(r)\leq(1+\eps_n)\cdot\delta_\kappa(r)+\eps_n\] on 
$(0,r_0)$ for some $r_0\in(0,2D_\kappa]$ and $\eps_n\to 0$. Then any ultralimit
$X_\om$ is locally CAT($\kappa$). More precisely, every closed ball in $X_\om$ of radius at most $\frac{r_0}{4}$ is
convex and CAT($\kappa$).
\ethm

\proof
By Proposition~\ref{prop_thickening}, 
there exists a sequence $(Y_n)$ of  complete length spaces and a sequence $(\beta_n)$
tending to zero such that  $Y_n$ is
a universal $\beta_n$-thickening of $X_n$ with the following properties.
$Y_n$
is $L$-Lipschitz $1$-connected up to some scale  for some universal $L\geq 1$. 
Moreover, there exist $C\geq 1$ and $r_1>0$ such that
each $Y_n$ admits a $(C,r_1)$-isoperimetric inequality  and  
$$\delta_{Y_n}(r)\leq(1+\beta_n)\cdot \delta_\kappa(r)+ \beta_n r^2$$
holds for all $r\in(\beta_n, r_0)$. Note that since $\beta_n\to 0$ we have
$\wlim Y_n= X_\om$. Let $\triangle\subset X_\om$ be a Jordan triangle of perimeter strictly less than $r_0$. By Theorem~\ref{thm_plateau},
we find an ultralimit $\hat Y_\om$ of a subsequence of $(Y_n)$ which 
contains an isometric copy of $\triangle$ and has the following property.
There exists a continuous map $v:\bar D\to \hat Y_\om$
which is a solution of the Plateau problem for $\triangle$ in 
$\hat Y_\om$  and satisfies
\[\a(v|_\Om)\leq \delta_\kappa(\ell(v|_{\D\Om}))\]
for every Jordan domain $\Om\subset D$ such that  
$\ell(v|_{\D\Om})<r_0$.
By Theorem~\ref{thm_int_top}, the associated intrinsic minimal disc $Z_{v}$ is homeomorphic to $\bar D$
and for every Jordan domain $\Om\subset Z_v$ with $\ell(\D\Om)<2D_\kappa$
holds
\[\hm^2(\Om)\leq\delta_\kappa(\ell(\D\Om)).\]
It follows from  the proof of \cite[Theorem~1.4]{LWcan} that the Dehn function of $Z_{v}$ is bounded above by $\delta_\kappa$ on $(0,2D_\kappa)$.
Thus, by \cite[Theorem~1.4]{LWcurv},  $Z_{v}$ is a CAT($\kappa$) space.
Let $u=\bar u\circ P$ be the induced factorization, see Section~\ref{sec_int_min}.
By Reshetnyak's majorization theorem,  $\D Z_{v}$ admits a $\kappa$-majorization $\varphi:C\to Z_{v}$
and then $\bar u\circ\varphi$ provides a $\kappa$-majorization for $\triangle$.
In particular, $\triangle$ satisfies the CAT($\kappa$) triangle comparison and
Lemma~\ref{lem_loc_cat}
implies the claim.
\qed

\bthm\label{thm_tree}
Let  $\delta:(0,r_0)\to \R$ be a continuous non-decreasing function
with 
\[\limsup_{r\to 0}\frac{\delta(r)}{r^2}<\frac{1}{4\pi}.\]
Suppose that  $(X_n)$ is a sequence of complete length spaces
such that the Riemannian Dehn functions satisfy 
\[\delta^{\mu^i}_{X_n}(r)\leq\delta(r)+\eps_n\] on 
$(0,r_0)$ for some sequence $\eps_n\to 0$. Then any ultralimit
$X_\om$ is 1-dimensional. More precisely, there exists $\tilde r>0$ depending only on 
the function $\delta$ such that every closed ball in $X_\om$ of radius at most $\tilde r$ is
convex and a tree.
\ethm

\proof
Let $X_\om$ be an ultralimit of the sequence $(X_n)$.
By assumption, for every $\ka\leq 0$ there exists $\tilde r_\ka>0$
such that $\delta(r)\leq\delta_\ka(r)$ holds for all $r\in(0,\tilde r_\ka)$.
From Theorem~\ref{thm_stable} we conclude that every closed ball of radius at most
$\frac{\tilde r_\ka}{4}$ in $X_\om$ is convex and CAT($\ka$). We set $\tilde r:=\frac{\tilde r_0}{4}$. Let $B\subset X_\om$ be a closed ball of radius at most $\tilde r$. Then $B$ is CAT(0) and therefore contractible. Moreover, $B$
is locally CAT($\ka$) for every $\ka<0$. We conclude from the Cartan-Hadamard theorem (\cite[Theorem~9.65]{AKP_found}) that $B$ itself is globally CAT($\ka$)
for every $\ka<0$. This completes the proof.
\qed
\medskip

Now we obtain Theorem~\ref{thm_asymp_cones}:
\bcor
Let $X$ be a complete length space  such that
\[\limsup\limits_{r\to\infty}\frac{\delta_X^{\mu^i}(r)}{r^2}\leq\frac{1}{4\pi}.\]
Then every asymptotic cone of $X$ is a CAT(0) space. Moreover, if the inequality is strict, then every asymptotic cone of $X$ is a tree. In particular, in this case, $X$
is Gromov hyperbolic.
\ecor

\proof
Set $C=\limsup\limits_{r\to\infty}\frac{\delta_X^{\mu^i}(r)}{r^2}$.
For every positive sequence $(\la_n)$ with $\la_n\to 0$ we find another positive 
sequence $(\eps_n)$ tending to zero such that 
\[\delta_{X_n}^{\mu^i}(r)\leq (1+\eps_n)C\cdot r^2+\eps_n \]
holds for all $r\geq 0$ and all $n\in\N$, where $X_n$ denotes the metric space
$(X,\la_n\cdot d)$. Thus Theorem~\ref{thm_stable} implies that any ultralimit of the sequence $(X_n)$ is CAT(0).
The additional statement in case $C<\frac{1}{4\pi}$ follows from Theorem~\ref{thm_tree}.
\qed
\medskip

It only remains to prove Theorem~\ref{thm_fine}.
\bthm\label{thm_asym_npc}
Let $\delta:(0,r_0)\to \R$ be a continuous non-decreasing function with
\[\limsup\limits_{r\to 0}\frac{\delta(r)-\frac{r^2}{4\pi}}{r^4}\leq 0.\] 
Suppose that  $(X_n)$ is a sequence of complete length spaces
such that the Riemannian Dehn functions satisfy 
\[\delta^{\mu^i}_{X_n}(r)\leq\delta(r)+\eps_n\] on 
$(0,r_0)$ for some sequence $\eps_n\to 0$. Then any ultralimit
$X_\om$ is locally CAT(0). More precisely, there exists $\tilde r>0$ depending only on 
the function $\delta$ such that every closed ball in $X_\om$ of radius at most $\tilde r$ is
convex and CAT(0).
\ethm

Note that the assumption on $\delta(r)$ cannot be relaxed to 
$\limsup\limits_{r\to 0}\frac{\delta(r)}{r^2}\leq\frac{1}{4\pi}$.
The latter condition merely ensures property (ET) and is satisfied by every Riemannian manifold.

\proof
Let $X_\om$ be an ultralimit of the sequence $(X_n)$.
Note that for $\ka>0$ we have
\[\delta_\ka(r)=\frac{1}{4\pi}\cdot r^2+\frac{\ka}{64\pi^3}\cdot r^4+o(r^5).\]
By assumption, for every $\ka> 0$ there exists $\tilde r_\ka>0$
such that $\delta(r)\leq\delta_\ka(r)$ holds for all $r\in(0,\tilde r_\ka)$.
From Theorem~\ref{thm_stable} we conclude that every closed ball of radius at most
$\frac{\tilde r_\ka}{4}$ in $X_\om$ is convex and CAT($\ka$). We set 
$\tilde r:=\frac{\tilde r_1}{4}$. Let $B\subset X_\om$ be a closed ball of radius at most $\tilde r$. Then $B$ is CAT(1) and  contractible. We claim that $B$
is CAT(0). By Cartan-Hadamard, it is enough to show that $B$ is locally CAT(0).
We will show the stronger statement that
every rectifiable Jordan curve $\Ga\subset B$ of length less than $2\pi$
admits a majorization.
Since $B$ is CAT(1) we find a continuous harmonic solution $u:\bar D\to B$ to the Plateau problem for $\Ga$
in $B$ and the associated intrinsic minimal disc $Z_u$ is a CAT(1) disc \cite[Appendix~A]{LWcurv}. In particular,  we have $\a(u)\leq 2\pi$.

We will finish the proof by showing that $Z_u$ is CAT(0).
Then the induced 1-Lipschitz map $\bar u:Z_u\to B$ leads the desired  majorization
as in the proof of Theorem~\ref{thm_stable}.

To show that $Z_u$ is CAT(0) it is enough to prove that $Z_u$ is CAT($\ka$) for every
$\ka\in(0,1)$, see \cite[Proposition~9.7]{AKP_found}.
Let $\ka\in(0,1)$ and cover the compact image of $u$ by a finite number of open balls $B_i$ such that each $B_i$ is
CAT($\ka$). This induces an open cover of $Z_u$. By compactness, there exists $\rho_\ka>0$
such that each ball of radius at most $\rho_\ka$ in $Z_u$ maps to one of the balls $B_i$. Since $B_i$ is CAT($\ka$),
it follows that $Z_u$ is locally CAT($\ka$) \cite[Corollary~1.2]{LWa_curv}.
Because $\hm^2(Z_u)=\a(u)\leq 2\pi$, \cite[Proposition~4.4]{LS_imp} implies that
$Z_u$ is a CAT($\ka$) disc. This completes the proof.
\qed

\bibliographystyle{alpha}
\bibliography{iso}

\begin{thebibliography}{HKST15}

\bibitem[AB04]{AB04}
S.~B. Alexander and R.~L. Bishop.
\newblock Curvature bounds for warped products of metric spaces.
\newblock {\em Geom. Funct. Anal.}, 14(6):1143--1181, 2004.

\bibitem[AB16]{AB16}
S.~B. Alexander and R.~L. Bishop.
\newblock Warped products admitting a curvature bound.
\newblock {\em Adv. Math.}, 303:88--122, 2016.

\bibitem[AKP23]{AKP_found}
S.~Alexander, V.~Kapovitch, and A.~Petrunin.
\newblock Alexandrov geometry: foundations.
\newblock {\em arXiv:1903.08539}, 2023.

\bibitem[BBI01]{BBI}
D.~Burago, Y.~Burago, and S.~Ivanov.
\newblock {\em A course in metric geometry}, volume~33 of {\em Graduate Studies
  in Mathematics}.
\newblock American Mathematical Society, Providence, RI, 2001.

\bibitem[BI12]{BI_min}
D.~Burago and S.~Ivanov.
\newblock Minimality of planes in normed spaces.
\newblock {\em Geom. Funct. Anal.}, 22(3):627--638, 2012.

\bibitem[BR33]{BR_sub}
E.~F. Beckenbach and T.~Rad\'{o}.
\newblock Subharmonic functions and surfaces of negative curvature.
\newblock {\em Trans. Amer. Math. Soc.}, 35(3):662--674, 1933.

\bibitem[EG15]{Evans}
L.~C. Evans and R.~F. Gariepy.
\newblock {\em Measure theory and fine properties of functions}.
\newblock Textbooks in Mathematics. CRC Press, Boca Raton, FL, revised edition,
  2015.

\bibitem[FW20]{FW_morrey}
M.~Fitzi and S.~Wenger.
\newblock Morrey's {$\varepsilon$}-conformality lemma in metric spaces.
\newblock {\em Proc. Amer. Math. Soc.}, 148(10):4285--4298, 2020.

\bibitem[GW20]{GW_noncompact}
C.-Y. Guo and S.~Wenger.
\newblock Area minimizing discs in locally non-compact metric spaces.
\newblock {\em Comm. Anal. Geom.}, 28(1):89--112, 2020.

\bibitem[HKST15]{HKST15}
J.~Heinonen, P.~Koskela, N.~Shanmugalingam, and J.~T. Tyson.
\newblock {\em Sobolev spaces on metric measure spaces}, volume~27 of {\em New
  Mathematical Monographs}.
\newblock Cambridge University Press, Cambridge, 2015.
\newblock An approach based on upper gradients.

\bibitem[Iva08]{I_vol}
S.~V. Ivanov.
\newblock Volumes and areas of {L}ipschitz metrics.
\newblock {\em Algebra i Analiz}, 20(3):74--111, 2008.

\bibitem[Kar07]{Kar_sobolev}
M.~B. Karmanova.
\newblock Area and co-area formulas for mappings of the {S}obolev classes with
  values in a metric space.
\newblock {\em Sibirsk. Mat. Zh.}, 48(4):778--788, 2007.

\bibitem[KS93a]{KS}
N.~Korevaar and R.~Schoen.
\newblock Sobolev spaces and harmonic maps for metric space targets.
\newblock {\em Comm. Anal. Geom.}, 1(3-4):561--659, 1993.

\bibitem[KS93b]{KS93}
N.~J. Korevaar and R.~M. Schoen.
\newblock Sobolev spaces and harmonic maps for metric space targets.
\newblock {\em Comm. Anal. Geom.}, 1(3-4):561--659, 1993.

\bibitem[LS19]{LS_conf}
A.~Lytchak and S.~Stadler.
\newblock Conformal deformations of {$\rm CAT(0)$} spaces.
\newblock {\em Math. Ann.}, 373(1-2):155--163, 2019.

\bibitem[LS20]{LS_imp}
A.~Lytchak and S.~Stadler.
\newblock Improvements of upper curvature bounds.
\newblock {\em Trans. Amer. Math. Soc.}, 373(10):7153--7166, 2020.

\bibitem[LS23]{LS_curv}
A.~Lytchak and S.~Stadler.
\newblock Curvature bounds of subsets in dimension two.
\newblock {\em J. Differential Geom., to appear}, 2023.

\bibitem[LW16]{LW_harm}
A.~Lytchak and S.~Wenger.
\newblock Regularity of harmonic discs in spaces with quadratic isoperimetric
  inequality.
\newblock {\em Calc. Var. Partial Differential Equations}, 55(4):Art. 98, 19,
  2016.

\bibitem[LW17a]{LW_plateau}
A.~Lytchak and S.~Wenger.
\newblock Area minimizing discs in metric spaces.
\newblock {\em Arch. Ration. Mech. Anal.}, 223(3):1123--1182, 2017.

\bibitem[LW17b]{LW_energy}
A.~Lytchak and S.~Wenger.
\newblock Energy and area minimizers in metric spaces.
\newblock {\em Adv. Calc. Var.}, 10(4):407--421, 2017.

\bibitem[LW18a]{LWint}
A.~Lytchak and S.~Wenger.
\newblock Intrinsic structure of minimal discs in metric spaces.
\newblock {\em Geom. Topol.}, 22(1):591--644, 2018.

\bibitem[LW18b]{LWcurv}
A.~Lytchak and S.~Wenger.
\newblock Isoperimetric characterization of upper curvature bounds.
\newblock {\em Acta Math.}, 221(1):159--202, 2018.

\bibitem[LW20]{LWcan}
A.~Lytchak and S.~Wenger.
\newblock Canonical parameterizations of metric disks.
\newblock {\em Duke Math. J.}, 169(4):761--797, 2020.

\bibitem[LW23]{LWa_curv}
A.~Lytchak and S.~Wagner.
\newblock Curvature bounds on length-minimizing discs.
\newblock {\em Preprint, arXiv:2308.14684}, 2023.

\bibitem[LWY20]{LWY_dehn}
A.~Lytchak, S.~Wenger, and R.~Young.
\newblock Dehn functions and {H}\"{o}lder extensions in asymptotic cones.
\newblock {\em J. Reine Angew. Math.}, 763:79--109, 2020.

\bibitem[Oss78]{O_iso}
Robert Osserman.
\newblock {The isoperimetric inequality}.
\newblock {\em Bulletin of the American Mathematical Society}, 84(6):1182 --
  1238, 1978.

\bibitem[PS19]{PS}
A.~Petrunin and S.~Stadler.
\newblock Metric-minimizing surfaces revisited.
\newblock {\em Geom. Topol.}, 23(6):3111--3139, 2019.

\bibitem[Res61]{R_iso}
Ju.~G. Reshetnyak.
\newblock An isoperimetric property of two-dimensional manifolds of curvature
  not greater than {$k$}.
\newblock {\em Vestnik Leningrad. Univ.}, 16(19):58--76, 1961.

\bibitem[Res68]{Resh-majorization}
Y.~Reshetnyak.
\newblock Non-expanding maps in a space of curvature no greater than $k$.
\newblock {\em Siberian Math. J.}, 9:918--927, 1968.

\bibitem[Res07]{Res97}
Yu.~G. Reshetnyak.
\newblock Sobolev-type classes of mappings with values in metric spaces.
\newblock In {\em The interaction of analysis and geometry}, volume 424 of {\em
  Contemp. Math.}, pages 209--226. Amer. Math. Soc., Providence, RI, 2007.

\bibitem[Ric21]{RR_2D}
Russell Ricks.
\newblock Closed subsets of a {${\rm CAT}(0)$} 2-complex are intrinsically
  {${\rm CAT}(0)$}.
\newblock {\em Algebr. Geom. Topol.}, 21(4):1723--1744, 2021.

\bibitem[Sta21]{St}
S.~Stadler.
\newblock The structure of minimal surfaces in {CAT}(0) spaces.
\newblock {\em J. Eur. Math. Soc. (JEMS)}, 23(11):3521--3554, 2021.

\bibitem[Tuk80]{Tuk80}
P.~Tukia.
\newblock The planar {S}ch\"{o}nflies theorem for {L}ipschitz maps.
\newblock {\em Ann. Acad. Sci. Fenn. Ser. A I Math.}, 5(1):49--72, 1980.

\bibitem[Wen08]{W_gromov}
S.~Wenger.
\newblock Gromov hyperbolic spaces and the sharp isoperimetric constant.
\newblock {\em Invent. Math.}, 171(1):227--255, 2008.

\bibitem[Wen11]{W_asym}
S.~Wenger.
\newblock The asymptotic rank of metric spaces.
\newblock {\em Comment. Math. Helv.}, 86(2):247--275, 2011.

\bibitem[Wen19]{W_dehn}
S.~Wenger.
\newblock Spaces with almost {E}uclidean {D}ehn function.
\newblock {\em Math. Ann.}, 373(3-4):1177--1210, 2019.

\end{thebibliography}

\end{document}